\documentclass[11pt,a4paper]{article}
\usepackage{amsmath,amsthm}
\usepackage{times,graphicx}

\usepackage{graphicx}
\usepackage{amssymb}
\usepackage{amsmath}
\usepackage{bm}
\usepackage{bbm}
\usepackage{setspace}
\usepackage{moreverb}

\usepackage[mathcal]{eucal}
\usepackage{mathrsfs} %per usare il carattere mathscr
\usepackage{verbatim}

\providecommand{\N}{\mathcal{N}}

\renewcommand{\Re}{\mathbb{R}}

\providecommand{\bx}{{\bf x}}

\providecommand{\bbeta}{{\bm \beta}}
\providecommand{\bdelta}{{\bm \delta}}

\def\d{{\rm d}}

\def\E{{\rm E}}

\def\N{\mbox{N}}

\def\d{{\rm d}}

\def\section{\bigskip}
\def\subsection{\medskip}

\newtheorem*{theorem*}{Theorem}

\begin{document}

\title{A General Framework for Updating Belief Distributions}
\author{P.G. Bissiri, C.C. Holmes \& S.G. Walker
\footnote{
Pier Giovanni Bissiri is Research Associate, Universita degli Studi di Milano-Bicocca,
Italy, (email: pier.bissiri@unimib.it);
Chris Holmes is Professor of Statistics, Department of Statistics, University of Oxford,
Oxford, U. K.  (email: c.holmes@stats.ox.ac.uk);
Stephen G. Walker is Professor, Department of Mathematics, University of Texas at Austin, U. S. A.  (email: s.g.walker@math.utexas.edu). Holmes is supported by the Oxford-Man Institute, Oxford, and the Medical Research Council, UK.
}}

\date{}
\maketitle

\vspace{0.1in}
\abstract{We propose a framework for general Bayesian inference. We argue that a valid update of a prior belief distribution to a posterior can be made for parameters which are connected to observations through a loss function rather than the traditional likelihood function, which is recovered under the special case of using self information loss. 

Modern application areas make it is increasingly challenging for Bayesians to attempt to model the true data generating mechanism. Moreover, when the object of interest
 is low dimensional, such as a mean or median, it is cumbersome to have to achieve this via a complete model for the whole data distribution. More importantly, there are settings where the parameter of interest does not directly index a family of density functions and thus the Bayesian approach to learning about such parameters is currently regarded as problematic. 

Our proposed framework uses loss-functions to connect information in the data to functionals of interest.  The updating of beliefs then follows from a decision theoretic approach involving cumulative loss functions. Importantly, the procedure coincides with Bayesian updating when a true likelihood is known, yet provides coherent subjective inference in much more general settings. Connections to other inference frameworks are highlighted. }

\vspace{0.1in} \noindent Keywords: Bayesian updating; PAC-Bayes, Decision theory; Generalized estimating equations; Gibbs posteriors; Information; Loss function; Maximum entropy; Self--information loss function.

\vspace{0.3in}
\noindent {\bf 1. Introduction.} Data sets are  increasing in size and modelling environments are becoming more complex. This presents opportunities for Bayesian statistics but also major challenges, perhaps the greatest of which is the requirement to define the true sampling distribution, or likelihood, for the data generator $f_0(x)$, regardless of the study objective. Even if the task is inference for a low-dimensional statistic of the population, Bayesian analysis is required to model the complete data distribution and, moreover, assume that the model is ``true''. 

In this paper we present a coherent procedure for general Bayesian inference which is based on the updating of a prior belief distribution to a posterior when the parameter of interest is connected to observations via a loss function. Briefly here, and in the simplest scenario, suppose interest is in the $\theta$ minimizing the expected loss
\begin{equation}
L(\theta)=\int l(\theta,x)\,\d F_0(x),\label{minloss}
\end{equation}
for some loss function $l(\theta,x)$, e.g. $l(\theta,x)=|\theta-x|$ for estimating a median, where $F_0(x)$ is the unknown distribution function from which i.i.d. observations arise. If $\pi(\theta)$ represents prior beliefs about this $\theta$, and $x$ is observed from $F_0$, then we argue that a valid and coherent update of $\pi(\cdot)$ is to the posterior $\pi(\cdot|x)$, where
\begin{equation}
\pi(\theta|x)\propto \exp\{-l(\theta,x)\}\,\pi(\theta).\label{bayest}
\end{equation}
The argument for this is given later in the paper, and to some extent relies on the idea that such an update {\sl must} exist. For we have a well defined parameter of interest $\theta$, an initial belief distribution about the location of the parameter, $\pi(\theta)$, and gain further information about $\theta$ via $x$ coming from $F_0(x)$. Updating beliefs is mandatory and it is obvious for some function $\psi$ we must have
$$\pi(\theta|x)=\psi(l(\theta,x),\pi(\theta)).$$
That the form for $\psi$ is (\ref{bayest}) is detailed later and a coherence property plays a key role:
\begin{equation}
\psi\bigg(l(\theta,x_2),\psi(l(\theta,x_1),\pi(\theta))\bigg)=\psi\big(l(\theta,x_1)+l(\theta,x_2),\pi(\theta)\big).\label{coher}
\end{equation}
This ensures we end up with $\pi(\theta|x_1,x_2)$ as the same object whether we update with  $(x_1,x_2)$ together  or one after the other.

A special case is when it is known that $F_0(x)=F(x;\theta_0)$ for some parametric family of distributions $F(\cdot;\theta)$, with corresponding density function $f(\cdot;\theta)$, and $l(\theta,x)=-\log f(x;\theta)$.
For minimizing $L(\theta)$ here yields $\theta_0$ and the update (\ref{bayest}) is the usual Bayesian update. It is important to note that our general update using loss functions should not to be seen as an approximation to anything; rather, it is targeting the parameter of interest, employing the necessary loss function with a valid update of beliefs.   

% For unless a probability model is ``correct", i.e. $F_0(x)=F(x;\theta_0)$, and this family is known, there is actually no formal connection between any observation $x$ and any parameter value of some chosen model. On the other hand, it is always possible to connect observation and parameter, even one not indexing a probability model, via a loss function. In fact, the parameter of interest is typcially defined by the loss function; see (\ref{minloss}). 

Classical inference based on the likelihood function can be regarded as using  the ``negative log likelihood function"  as a loss function; for example, in the case of independent and identically distributed observations, we can regard
$$L(\theta;x_1,\ldots,x_n)=-\sum_{i=1}^n \log f(x_i|\theta)$$
as a loss function connecting data $(x_i)$ with a parameter $\theta$ indexing the family of density functions $f(x|\theta)$. And in this setting, one does not even need to assume the correctness of the model; one is merely expressing interest in the parameter $\theta_0$ minimizing
$$-\int \log f(x;\theta)\,\d F_0(x)$$
which is the parameter minimizing the Kullback-Leibler divergence between the unknown $f_0(\cdot)$ and the family $f(\cdot;\theta)$.

\vspace{0.2in}
\noindent
{\bf 1.1 The idea.}  Here we provide further elaboration on the outline of the idea given previously. Let $\theta$ denote a parameter or functional of interest, for example the mean or median of a population $F_0(x)$, and let $x$ denote an  observation from $ F_0(x)$, with $F_0$ unknown. We are interested in a formal way to update prior beliefs $\pi(\theta)$ to posterior beliefs $\pi(\theta | x)$ given $x$. 

Bayesian inference proceeds through knowledge of a complete, true, model for $f_0(x)$. This is often parameterised via a sampling distribution $f(x; \theta)$ and a prior $\pi(\theta)$, and define the marginal likelihood
$$m(x) = \int f(x;\theta) \pi(\d \theta).$$
Then, see for example Bernardo and Smith (1994),
inference  for $\theta$ occurs via Bayes theorem
$$ \pi(\theta | x) = f(x ; \theta) \pi(\theta) / m(x) .$$
However, the statement ``inference for $\theta$" is meaningless unless the true parametric family $f(\cdot;\theta)$ is known.
In this case, following the Savage axioms (Savage, 1954), the Bayesian update can be shown to be the rational way to proceed. However, $f_0(x)$ may be unknown,  and even if $f(\cdot;\theta)$ is correct, $\theta$ might be ultra high dimensional mainly made up of nuisance parameters relative to a low dimensional subset of the parameters of interest. Taken all together, these points can make the Bayesian approach cumbersome.

We are interested in the rational updating of beliefs under more general and less stringent conditions. To do so we make use of loss functions to connect information in data to parameters of interest.  Informally for now, we write such loss functions as $l(\theta,x)$, and we will discuss specific types later in the paper.  We shall consider the reporting of subjective beliefs, $\pi(\theta | x)$, as an action made under uncertainty and use decision theory to guide the optimal action.  See, for example, Hirshleifer and Riley (1992).

To outline the theory, let $\nu$ denote a probability measure on the space of $\theta$. We shall construct a loss function to select an optimal posterior distribution $\widehat{\nu}(\theta)$\footnote{We use $\widehat{\nu}$ to denote optimality rather than an approximation or estimate} given a prior $\pi(\theta)$ and data $x$. To achieve this we construct a loss-function $L(\nu ; \pi, x)$ on the space of probability measures on $\theta$ space, and then present 
$$
\widehat{\nu} = \arg \min_{\nu} L(\nu; \pi, x),
$$
as the representation of beliefs about the unknown value of $\theta$ given the prior information, represented via the belief distribution $\pi$, and data $x$. 
As it is widely assumed data $x$ is an independent piece of information to that which gave rise to the prior, it is appropriate to consider an additive, or cumulative, loss function of the form
\begin{equation}\label{eq:log-form}
L(\nu; \pi, x)=h_1(\nu,x)+h_2(\nu,\pi),
\end{equation}
 where $h_1$ and $h_2$ are themselves loss functions on probability measures,  representing fidelity-to-data and fidelity-to-prior, respectively. See, for example, Berger (1993) for more about ideas on uses of loss functions within decision theory.

The question is whether we can claim a probability measure selected as the solution to a decision problem; i.e. minimizing a loss function, can be viewed as representing beliefs about a parameter. To answer this, given the aim (\ref{minloss}), we would clearly prefer probability measure $\nu_1$ to $\nu_2$ as representing beliefs if
\begin{equation}
\int\int l(\theta,x)\,\d F_0(x)\,\nu_1(\d\theta)\leq \int\int l(\theta,x)\,\d F_0(x)\,\nu_2(\d\theta).\label{prefer}
\end{equation}
Indeed, it would be incoherent to select $\nu_2$ rather than  $\nu_1$ when (\ref{prefer}) holds. 
Thus the answer is affirmative.  Though we are not minimizing or comparing (\ref{prefer}), since we do not have $F_0$, we can substitute the expression
\begin{equation}
L(\nu;F_0) =  \int\int l(\theta,x)\d F_0(x)\,\nu(\d\theta)\label{asloss} 
\end{equation}
with the Bayesian finite sample expression of the form (\ref{eq:log-form}). We now discuss the choices of $h_1$ and $h_2$ which give (\ref{eq:log-form}) as a  Bayesian finite sample version of (\ref{asloss}).

Under this approach the analyst needs to specify $h_1$ and $h_2$ in such a way that they proceed in an optimal, rational, and coherent manner. 
% We can deal immediately with the loss function $h_2(\nu,\pi)$. 
Somewhat remarkably, as proved in the Supplementary Material, for coherent inference (\ref{coher}), $h_2$ must be the Kullback--Leibler divergence, Kullback and Leibler (1951), and given by
$$h_2(\nu,\pi)=d_{KL}(\nu,\pi)=\int \nu(\d\theta)\,\log\{\nu(\d\theta)/\pi(\d\theta)\}.$$
%whenever $\nu$ is absolutely continuous with respect to $\pi$, which we will always assume.
%Since, in reality, the information $I$ is pertaining to $\theta$, it is convenient and realistic to construct a loss function $l(\theta,I)$ where, as is usual with loss functions, $l(\theta,I)$ denotes the loss in taking $\theta$ as the optimal value when presented with information $I$. 
Regarding $h_1$, since $\nu(\theta)$ is a probability measure representing beliefs about $\theta$, the only choice here is to take the loss-to-data $h_1(\nu,x)$ on the probability measure as the {\sl expected} loss (see von Neumann and Morgenstern, 1944) of  $l(\theta,x)$; that is
$$h_1(\nu,x)=\int l(\theta, x)\, \nu(\d\theta),$$
with the particular types of the loss-function on the parameter on interest $l(\theta,x)$ to be discussed later. 
% In general there the form of $l(\theta, x)$ will be problem specific as discussed in Section 3.
% That there is no ``right" loss function should not be seen as a concern for the choice of $h_1$. However, in this case, there is no alternative choice to $h_1$ which could even be proposed. 

Substituting in $h_1$ and $h_2$, the cumulative loss function is then given by
\begin{equation}\label{f:cum_loss}
L(\nu;\pi,x)=\int l(\theta,x)\,\nu(\d\theta)+d_{KL}(\nu,\pi).\end{equation}
This then; i.e. (\ref{f:cum_loss}), is our finite sample version of (\ref{asloss}), and note that (\ref{f:cum_loss}) becomes (\ref{asloss}) as $n\rightarrow\infty$.  The solution to (\ref{f:cum_loss}) provides the $\widehat{\nu}$ which the statistician believes best minimizes (\ref{asloss}). This is, according to our approach, done by using the empirical distribution function as a substitute for $F_0$ and using a penalty term which prevents the answer from being too far from the prior in a Kullback-Leibler sense; the Kullback-Leibler appearing here for the necessary coherence property of the answer. Of interest, as discussed later on, is the PAC Bayes solution to the problem (Langford, 2005) that finds an approximation which minimizes an upper bound for (\ref{asloss}). 

Surprisingly, but quite easy to show, the minimizer of $L(\nu; \pi, x)$ is given by
\begin{eqnarray} \label{eq:bayes}
\widehat{\nu}(\theta) & = & \arg \min_{\nu} L(\nu; \pi, x)  \nonumber \\
& = & \frac{\exp\{-l(\theta,x)\}\pi(\theta)}{\int \exp\{-l(\theta,x)\}\pi(\d\theta) }.
\end{eqnarray}
This can be seen by observing that
$$\int l(\theta,x)\,\nu(\d\theta)+d_{KL}(\nu,\pi)=\int \nu(\d\theta)\log\left\{\frac{\nu(\theta)}{\exp(-l(\theta,x))\pi(\theta)}\right\}.$$
So (\ref{eq:bayes}) has the form of a Bayesian update.  As is usual in decision problems involving the use of loss functions, it is incumbent on the decision maker to ensure solutions exist. So $l(\theta,x)$ needs to be constructed such that
$$0<\int \exp\{-l(\theta,x)\}\pi(\d\theta)<+\infty.$$
Whereas the Bayesian approach requires the construction of a probability model for all possible outcomes conditional on all unknown states of nature, the  approach here  requires the construction of loss functions given the outcomes for only the parameter of interest. This allows the decision maker to concentrate on modeling only those quantities that are important to the task at hand.

\vspace{0.2in}
\noindent
{\bf 1.2 Connections with related work.}
There is a large literature on procedures for robustly estimating a parameter of interest by minimizing the cumulative loss
\begin{equation}
L(\theta;x)=\sum_{i=1}^n l(\theta,x_i).\label{fcumloss}
\end{equation}
This is clearly the finite sample version of
$$L(\theta)=\int l(\theta,x)\,\d F_0(x).$$
Our claim is that (\ref{f:cum_loss}) is the Bayesian version of (\ref{fcumloss}), where interest is on probability measues on $\theta$ space rather than single states $\theta$.

H\"uber (2009) provides examples of (\ref{fcumloss}), where we note that the primary aim is not modeling the data but rather estimating a statistic. This is an advantage when a probability model for the data is too hard to formulate.
We are presenting a Bayesian extension of this idea. Since we are interested in a belief distribution for $\theta$ given data, and have further information provided by $\pi$,  we claim the appropriate Bayesian version is given by (\ref{eq:bayes}).

Some of the ideas presented in the paper have been considered in  a less general setting by Zhang (2006a, 2006b) and  
Jiang and Tanner (2008). In Zhang (2006a) an estimation procedure,  named Information Risk Minimization, also known as a Gibbs posterior, which has the same form as \eqref{eq:bayes}, is described in Section IV of his paper. This is our procedure when data is regarded as stochastic. Zhang then  concentrates on the  properties of the Gibbs posterior. 

Further theoretical work is done in Zhang (2006b). In Jiang and Tanner (2008) a Gibbs posterior is studied in comparison with a true Bayesian posterior where the model is assumed to be misspecified. The claim is that posterior performance of a Bayesian model can be unreliable when misspecified, whereas a Gibbs posterior which targets points of interest can have better performance. The comparison involves variable selection for high-dimensional classification problems involving a logit model.

We build on the work of Zhang (2006a, 2006b) and  
Jiang and Tanner (2008) in a number of important directions. The first is that we show that the Gibbs posterior is the only coherent, decision theoretic approach for inference and statistical applications under misspecification. We provide a principled approach to scale the relative information in the data to information in the prior (see Section 3); that is left as an arbitrary free parameter in Zhang (2006a, 2006b) and  
Jiang and Tanner (2008). 

Bissiri and Walker (2010) use (\ref{f:cum_loss}) with Bernoulli observations and find sufficient conditions on $l(\theta,x)$ for the sequence of posteriors, based on (\ref{eq:bayes}), to be consistent.
This result for consistency is extended to more general i.i.d. observations in Bissiri and Walker (2012a). 
In Bissiri and Walker (2012b) it is shown that starting from the class of $g$-divergences, for a coherent sequence of updates, which will be explained later in the paper, we need the Kullback-Leibler divergence as the loss between prior $\pi$ and $\nu$. In this paper, the Theorem presented in the Appendix, we provide a simpler and more intuitive proof to this result by assuming a straightforward and realisitic assumption about the class of $g$-divergences.

A similar construct to $L(\nu;\pi, x)$ is provided by Zellner (1988), who presents what is essentially a loss function for the posterior distribution using ideas of information processing from prior to posterior.  The motivation is different and relies on notions of information present in log probabilities and log likelihoods, which may not be compatible as noted by J.M. Bernardo in the discussion of Zellner's paper.
Furthermore, our derivation of the loss function allows a broader interpretation of the elements, which does not require the existence of a probability distribution for the observation.

Concerns that the specification of a complete model for the data generating distribution is unachievable date back to de Finetti (1937) and the notion of ``prevision''. In his work de Finetti considers conditional expectation as the fundamental primitive, 
or statistic, of interest on which prior beliefs are expressed and updated. Recently other researchers have further developed this approach under the field of Bayesian linear statistics, see Goldstein and Wooff (2007).

There has been increasing awareness of the restrictive assumptions that formal Bayesian analysis entails. Royall and Tsou (2003) describe procedures for adjusting likelihood functions when the model is misspecified. More recently, Doucet and Shepherd (2012), and Muller (2012) consider formal approaches to pseudo-Bayesian methods using sandwich estimators to update subjective beliefs, motivated by robustness to model misspecification, see also Hoff and Wakefield (2013). Ribatet et al. (2009) consider pseudo-Bayesian approaches with composite likelihoods. More generally there is increasing recognition that formal Bayesian analysis can be restrictive for example through computational issues, such as arise in the area of Approximate Bayesian Computation (see for example Marin {\em{et al}} (2012)).

Several authors have considered issues with Bayesian updating by using proxy models, $f(x; \theta)$, for example, see Key et al. (1999), when $(x_i)$ is known not to arise from $f(x;\theta)$ for any value of $\theta$. That is, there is no $\theta$ conditional on which $x$ is from  $f(x;\theta)$. This is referred to as the M--open case in Bernardo and Smith (1994). One suggested solution is to use methods based on approximations and Key et al. (1999) describe one such idea using a cross--validation approach. While this may be a pragmatic it does have some shortcomings. Most serious is that there is little back--up theory and this has repercussions in that the update suffers from a lack of coherence

Another approach is to ignore the problem. That is, assume the observations are coming from $f(x;\theta)$ even though it is known they are not. According to Goldstein (1981), ``there is no obvious meaning for Bayesian analysis in this case". The disaster of making horribly wrong inference can be protected to some extent by model selection; that is, postulating a number of models for $f_0(x)$, say $f_j(x;\theta_j)$, with corresponding priors $\pi_j(\theta_j)$, and model probabilities $(p_j)$, for $j=1,\ldots,M$. But as Key et al. (1999) point out, how does one construct $\pi_j(\theta_j)$ and $p_j$ when one knows none of the postulated models are correct. So the Bayesian update breaks down in that nothing has any interpretation.

Finally, and we acknowledge the contribution of the reviewers for pointing this out, we discuss connections with PAC Bayes; see, for example,  Shawe-Taylor and Williamson (1997), Langford (2005), Alquier (2008), McAllester (1998), and also the approach in Cesa-Bianchi and Lugosi (2006). PAC Bayes is an interesting emerging field in machine learning concerned with techniques for bounding the generalisation error (empirical risk) of a Bayesian model. The motivation behind PAC Bayes is to find an upper bound for the empirical risk of a probability measure $\nu$ on a model $L(\nu;F_0)$ in  (\ref{asloss}), termed generalisation error in the PAC Bayes literature. Given observation $x$ and prior $\pi$,  the upper bound will be written as $U(\nu;x,\pi)$; i.e., for all $\nu$ it is that
$$L(\nu;F_0)\leq U(\nu;x,\pi).$$
See Catoni (2003) where the form of $U$ is provided. Then it can be shown that an upper bound
$U(\nu;x,\pi)$ is provided by (\ref{eq:bayes}). The PAC Bayes approach is complementary to our work. The motivation and construction is very different. We are interested in a framework for the rational updating of beliefs, rather than seeking bounds on the empirical risk of a probability measure on models. The minimizer of an upper bound is interesting but does not justify using $\widehat{\nu}$ as an update of a belief distribution for Bayesian style inference, and hence whether $\widehat{\nu}$ form a coherent sequence of belief distributions is not discussed in the PAC Bayes formulation of $U$. Whereas the requirement of coherence is central to our formulation. Moreover the scaling of the loss-to-data $h_1$ to the loss-to-prior $h_2$ enters as a constant in the margin of the error bound in PAC Bayes whereas here it has explicit meaning in the relative weight of information provided by the two sources, prior and data (see Section 3).
% That the sequence is coherent is more by chance than design. 

This said, there are clear synergies and the operational characteristics of PAC Bayes are similar, they must be since we gather the same answer. However, the motivation is different and we sincerely believe our derivation is providing full support for the use of the posterior. Moreover, as we will see later, this derivation provides insights into the necessary calibration of loss functions $h_1$ and $h_2$.

\vspace{0.2in}
\noindent
{\bf 1.3 Layout of paper.}
The layout of the remainder of the paper is as follows. In Section 2 we discuss types of loss function. When the self-information loss function is used then the update is the traditional Bayes update. With other loss functions there is a calibration issue between the two styles of loss function used, i.e. the loss to the data and the loss to the prior.  This calibration problem is discussed and resolved in a number of ways in Section 3. In Section 4 we discuss forms of information other than the usual data arising from some unknown distribution function. This includes non-stochastic information and also partial information.
Section 5  provides some numerical illustrations including inference based on partial information  and a clustering problem.
Section 6 concludes with a discussion on a number of points.

%\vspace{0.2in}
%\noindent
%{\bf 2 Specification of loss functions.}
%Under our approach \eqref{eq:log-form}, loss-functions dictate how the analyst expresses beliefs on $\theta$ using loss functions to the data $h_1$ and loss to the prior $h_2$. In this section, we show how these loss functions should be chosen in order to ensure that the analyst remains coherent and honest.

\vspace{0.2in}

\noindent
{\bf 2. Types of loss function.} In this section we will consider the  form of $h_1$ in \eqref{eq:log-form} that connects information in the data to the value of the unknown $\theta$. We shall consider three broad situations. First, when the analyst believes they know the complete family of distributions from which the $(x_i)$ arose, the so called M--closed scenario. Second, when $f_0(x)$ is unknown but where a complete likelihood $f(x; \theta)$ is being used as a proxy model. Finally, when there is no  sampling distribution or proxy model for $x$ and the parameter of interest is connected to $x$ via a loss function $l(\theta,x)$.

\vspace{0.2in}
\noindent
{\bf 2.1 M--closed and self-information loss.}
When the analyst knows the family from which $(x_i)$ arose, the so-called M--closed view, then the Bayesian approach to learning is fully justified, well known and widely used as a statistical approach to inference; the book of Bernardo and Smith (1994) is comprehensive.  
%Here we recall the essence of it: A parameter of a density function $f(x;\theta)$, $\theta\in\Theta$, is unknown and beliefs about it are encapsulated via a prior distribution $\pi(\theta)$.  Once (conditionally) independent samples $(x_1,\ldots,x_n)$ are observed from the density function $f(x;\theta)$, the prior is updated to the posterior distribution via Bayes' Theorem; given by
%$$\pi(\d\theta|x_1,\ldots,x_n)=\frac{l_n(\theta)\,\pi(\d\theta) }{\int_\Theta l_n(\theta)\,\pi(\d\theta)},$$
%where $l_n(\theta)=\prod_{i=1}^n f(x_i;\theta)$ is the likelihood function.  The posterior then represents revised beliefs taking into account both the prior distribution and the observations. Mathematically, it is an application of Bayes' Theorem via the standard definition of conditional probability.
%
%So the Bayesian update works and is applicable in the case when the $(x_i)$ come from the density $f(x;\theta)$ for some $\theta\in\Theta$. 
% In Bernardo and Smith (1994) this is referred to as the M--closed view. 
To see how Bayes arises in our framework, we would need to construct a loss function for $l(\theta,x)$ with the knowledge that $x$ came from $f(x;\theta)$. It is well known that the appropriate and sole loss function in this case is the self--information, or logarithmic loss function, given by
$$l(\theta,x)=-\log f(x;\theta).$$
Indeed, the cumulative loss version of this is the log-likelihood function.

See Bernardo (1979) and Merhav and Feder (1998) for more on the self information loss function. This amounts to the use of proper scoring rules to ensure that the analyst remains honest in 
 expressing subjective beliefs when the parametric family $f(x; \theta)$ is known, under which our approach coincides  with the Bayesian updating rule.

 % However, there are different ideas behind our derivation of this rule, with different assumptions being made. Most crucially, we need the $(x_i)$ to provide independent pieces of information to maintain the credibility of the cumulative loss function.

\vspace{0.2in}
\noindent
{\bf 2.2 M--open and the use of proxy models.}
% As has been mentioned by many authors,, 
Issues with the Bayesian rule arise when the form of $f(x;\theta)$ is not known, for example, see Key et al. (1999). Equivalently, there is no $\theta$ conditional on which $x$ is from  $f(x;\theta)$; more bluntly, there is no connection between any $x$ and any $\theta$ via $f(x;\theta)$. This is referred to as the M--open case in Bernardo and Smith (1994). In many situations, the correct sampling density, $f_0(x)$, is unknown or unavailable or too complex to work with. 

One way to proceed is by considering $\theta_0$, the value of $\theta$ that minimizes the the Kullback--Leibler divergence between a proxy model $f(x; \theta)$ and the true density function $f_0(x)$; i.e. $\theta_0$ minimizes
$$d_{KL}(f_0(\cdot),f(\cdot;\theta))=\int f_0(x)\,\log \{f_0(x)/f(x;\theta)\}\,\d x.$$
%and we will let
%$$\delta=\int f_0(x)\,\log\{f_0(x)/f(x;\theta_0)\}\,\d x.$$
Then prior beliefs, $\pi(\theta)$, will be expressed on this unknown value. We show in Appendix B that it is then possible to learn about this $\theta_0$ since an infinite collection of $(x_i)$ yields $\theta_0$. Then we would wish the sequence of $\nu(\theta)$ to accumulate about $\theta_0$. The appropriate loss function in this case is still $l(\theta,x)=-\log f(x;\theta)$. For the standardized cumulative loss based on a sequence of observations is given by
$$-n^{-1}\sum_{i=1}^n \log f(x_i;\theta)\rightarrow -\int \log f(x;\theta)\,\d F_0(x)\quad\mbox{a.s.}\quad\mbox{for all}\,\,\theta,$$
which is minimized by $\theta_0$.

%
%When an approximate model $f(x;\theta)$ has been supposed, it is often prudent to consider a number of models, say $f_j(x;\theta_j)$ for $j=1,\ldots,M$, as we have mentioned previously. We can deal with this in a simple way. So let $\theta=(\theta_1,\ldots,\theta_M)$ and let $\pi(\theta)$ be the prior distribution for $\theta$ on $\Theta=\cup_{j=1}^M\Theta_j$. This would be constructed by considering beliefs about which $\theta_j$ from $f_j(\cdot;\theta_j)$ takes this family closest to $f_0(\cdot)$. The model $f(x;\theta)$ would then be given by
%$$f(x;\theta)=\sum_{j=1}^M p_j\,f_j(x;\theta_j)$$
%and the $(p_j)$ would now be the probabilities describing beliefs about which model provides the closest density to $f_0(\cdot)$. Hence, unlike the Bayesian approach to model selection in the M--open case, all the quantities to be specified have clear interpretation. We can recover the Bayesian update when we take, for each $i\in(1,\ldots,n)$,
%$$l(\theta,x_i)=-\log f(x_i;\theta).$$

So while the Bayesian approach has foundational issues to deal with whether the  M--open or M--closed view hold, for the approach here it is irrelevant. If one adopts $\theta_0$ as the parameter value taking the family closest to $f_0(\cdot)$ then one does not need to worry if one is M--open or M--closed, since if $f(\cdot;\theta)$ is the true family then obviously $\theta_0$ reverts to the true parameter value. This point is crucial, since for the Bayesian being in M--open or M--closed forces one to adopt different inference approaches, see Bernardo and Smith (1994).

Moreover our approach supports the use of the relevant partial information in the data for updating beliefs on the parameter of interest, such as in the proportional hazards model. This can be especially important when the data is high dimensional. Such updates have no formal justification from a Bayesian perspective.
% it may be that one simply does not know if one is in the M--open or M--closed view (though strictly speaking this puts you in M--open) and then one needs a framework in which the same approach is adopted and justified  regardless of which view is held. We have provided such a framework.

\vspace{0.2in}
\noindent
{\bf 2.3 Paremeter minimizing a loss function.} In the most general scenario the parameter of interest minimizes a loss function of the type (\ref{minloss}).  In the classical literature, 
this type of estimation problem is in the area of {\it Robust Statistics} and specific loss functions can be found in the literature, pertaining to $M$-estimation and estimating equations. 
See, for example, H\"uber (2009). 

An important class of loss functions is provided by the $M$ estimators for a location parameter, H\"uber (1964). So rather than using the loss function $-\log f(x_i;\theta)$, a $\rho(x_i;\theta)$ is used in an attempt to obtain robust estimation, rather than the traditional maximum likelihood estimator, which can be suspect if the model is incorrect. This idea has been generalized to the class of estimating equations, whereby the estimate of $\theta$ is obtained by minimizing
$$\sum_{i=1}^n \rho(x_i;\theta).$$
Our approach, which mirrors this classical robust procedure, would use the loss function
$$L(\nu;x_1,\ldots,x_n,\pi)=\int\sum_{i=1}^n \rho(x_i;\theta)\,\nu(\d\theta)+d_{KL}(\nu,\pi)$$
with solution provided by
$$\widehat{\nu}(\d\theta)\propto \exp\left\{-\sum_{i=1}^n \rho(x_i;\theta)\right\}\,\pi(\d\theta).$$
%For example, one possible application would be the Generalized Estimating Equations, see Liang and Zeger (1986). For the grouped observations $(x_{i1},\ldots,x_{in_i})$,
%$$\rho(x_i,\theta)=\half(x_i-\mu_i(\beta))'V_i(\phi,\alpha)^{-1}(x_i-\mu_i(\beta))$$
%where $\theta=(\beta,\phi,\alpha)$ and for some link function $g$, $g(\mu_{ij}(\beta))=x_{ij}'\beta$, and for some correlation matrix $R_i(\alpha)$ and diagonal matrix $A_i$, with $j$ entry
%given by $a(\mu_{ij}(\beta))$, with $a$ a specified variance function, $V_i=\phi A_i^{1/2}R_i(\alpha)A_i^{1/2}$, with $\phi$ a scale parameter.
% There is by now an abundance of literature on $M$-estimation, estimating equations and generalized estimating equations. Our point is that all such equations can be viewed as loss functions connecting independent units with parameters of interest. Hence, all fit within our framework and we would extend the loss function to include the prior $\pi$ and we obtain an explicit expression for $\nu(\d\theta)$. In cases when the parameter estimation is done via iterative methods, which is typically the case, Markov chain Monte Carlo methods would substitute for our sampling strategies from $\nu(\d\theta)$.

% In essence, this is the practical innovations of the framework we are proposing. We are claiming that any loss function of the type
% $$\sum_{i=1}^n \rho(x_i,\theta)$$
% can be extended to the Bayesian type updating mechanism. 
The $\theta_0$ of interest is implicitly assumed to be the limit of the sequence of minimizers of the cumulative losses. This would be the minimizer of $\int \rho(x;\theta)\,\d F_0(x)$ and hence the prior beliefs are being expressed about this unknown value. Then the loss function $l(\theta,x)=\rho(x;\theta)$ is ensuring the updates are indeed ``moving towards" $\theta_0$. To complete the picture, it would have been that the decision maker would be happy to make a decision given the minimizer of $\int \rho(x;\theta)\,\d F_0(x)$.

\vspace{0.2in}
\noindent
{\bf 3. Calibration of relative losses.} This section deals with the important aspect of specifying the relative information in the data to the information in the prior in general settings.  In the M--closed and M--open settings the use of the self-information loss $l(\theta, x)= -\log f(x;\theta)$ results in a fully specified form for \eqref{eq:bayes}. However, in the setting of Section 2.3 there is an issue about the scale of the loss function $h_1$ which is a consequence of the apparent arbitrariness in the weight of $l(\nu,x)$ relative to $l(\nu, \pi)$, in that we are free to multiply either by an arbitrary factor.  So, equivalently, we are interested in the loss function $w\,l(\theta,x)$ for some $w>0$. The question is how to select $w$, noting that $w$ controls the relative weight of loss-to-data to loss-to-prior. Of course, such an issue does not arise in the classical literature on estimation using such loss functions since there is no combining with different styles of loss functions. 
However, the calibration of different types of loss function is not a unique problem. It arises in many applied contexts; possibly the most well known be in health economics where losses pertaining to costs need to be balanced against losses pertaining to health benefits.

The most common ideas for assigning $w$ in Gibbs posteriors and PAC Bayes typically involve cross validation and subjective choices. As mentioned above, in PAC Bayes the weighting $w$ is a constant that enters into the margin of the error bound. Here we discuss some more ideas intended to help the analyst. We do not claim to be exhaustive in the approaches, or to be prescriptive in advocating one approach over an other. Our intention is to provide tools for elicitation of the relative loss-to-data to loss-to-prior.

\vspace{0.2in}
\noindent
{\bf 3.1 Annealing. } In the literature on Gibbs posteriors, the weighting parameter  is labelled as a ``temperature" and selected subjectively. 
There are clear connections here with the use of ``power priors'' (Ibrahim \& Chen, 2000) where
$$\nu(\d\theta)\propto \prod_{i=1}^n f(x_i;\theta)^w\,\pi(\d\theta).$$
Such an idea has also been discussed in Walker and Hjort (2001). It is evident what $w$ achieves; if $0<w<1$ then the loss-to-prior is given more prominence than in the Bayesian update and the data will be less influential. In the extreme case when $w=0$ we retain the prior throughout. On the other hand, when $w>1$ the loss $-\log f(x;\theta)$ is given more prominence than in the Bayesian update and in the extreme case when $w$ is very large the $\nu$ is accumulating about the maximum likelihood estimator for the model; that is $$\nu(\d\theta)\approx \delta_{\widehat{\theta}}(d\theta),$$ where $\widehat{\theta}$ maximizes $\prod_{i=1}^n f(x_i;\theta)$.

\vspace{0.2in}
\noindent
{\bf 3.2 Unit information loss. }
Here we discuss a procedure for default subjective assignment based on a prior evaluation of the expected value of $l(\theta,x)$. The idea originates from work 
in the specification of reference priors and ``objective Bayes'', see for example Kass and Wasserman (1996).

% While this can be chosen, an empirical calculation would use
% $$\int\int l(\theta,x)\,\pi(\d\theta)\,\d F_n(x)$$
% where $F_n$ is the empirical distribution function of the data.

To aid in the calibration of the loss functions and the selection of $w$ we can consider the following. To begin we need to ensure that both losses are non-negative for all $\theta$. Hence we write the loss function with an additional term $\log \pi(\widehat{\theta})$, which is a constant, and where $\widehat{\theta}$ maximizes $\pi(\theta)$, so that the cumulative loss becomes
$$L(\nu;x,\pi)=\int\left[w \,l(\theta,x)+\log\{\pi(\widehat{\theta})/\pi(\theta)\} \right]\nu(\d\theta)+\int\nu(\d\theta)\,\log \nu(\theta).$$
and we would additionally standardise $l(\theta, x)$ such that $\min_{\theta} l(\theta,x)=0$ for any $x$. If this is not the case then we replace $l(\theta,x)$ by $l(\theta,x)-l(\theta_x,x)$ where now $\theta_x$ minimizes $l(\theta,x)$.
Hence, we can regard
$$L(\theta;x,\pi)=w\,l(\theta,x)+\log\{\pi(\widehat{\theta})/\pi(\theta)\} $$
as a loss function for $\theta$ with information provided by $x$ and $\pi$.
So, assuming that $l(\theta,x)>0$, we want to calibrate the two loss functions given by
$$w\,l(\theta,x)\quad\mbox{and}\quad \log \{\pi(\widehat{\theta})/\pi(\theta)\}.$$

These are two loss functions for $\theta$ and to adhere with the notion that before we have any data, there is a single piece of information, we can calibrate the two losses by making the expected losses, taken over $\theta$ and $x$, to match. That is, whether someone takes a $\theta$ and is penalized by the loss $$\log\{\pi(\widehat{\theta})/\pi(\theta)\},$$ or
takes a $(\theta,x)$ and is penalized by the loss $wl(\theta,x)$, at the outset, the expected losses should match. They are confronted by two choices of loss with one piece of information and thus the losses can be calibrated by ensuring their expected losses coincide. The connection between expected information and expected loss can be found in Bernardo (1979).

Thus $w$ can be set by ensuring
$$w\E_{\theta, x} \left(l(\theta,x)\right) =\E_{\theta}\left(\log \{\pi(\widehat{\theta})/\pi(\theta)\}\right).$$
Here $\E$ is with respect to a joint belief in $x$ and $\theta$; say $m(x,\theta)$, the marginal for $\theta$ of which is $\pi(\theta)$.
So
$$w=\frac{\int\log \{\pi(\widehat{\theta})/\pi(\theta)\}\,\pi(\d\theta) }
{\int \int l(\theta,x)\,m(\d\theta,\d x) }.$$
This requires specifying a sampling distribution for $x$, which we've said is problematic,  so an empirical choice is then given by
$$w=\frac{\int\log \{\pi(\widehat{\theta})/\pi(\theta)\}\,\pi(\d\theta) }
{\int\int l(\theta,x)\,\pi(\d\theta)\,\d F_n(x)}.$$
where $F_n$ denotes the empirical distribution function.

Let us consider an example, where $l(\theta,x)=(\theta-x)^2$ with $\pi(\theta)=\N(\theta|0,\tau^2)$ with
$m(x|\theta)$ being any density with mean $\theta$ and variance $\sigma^2$.
Then we can evaluate
$$\int \log\{ \pi(\widehat{\theta})/\pi(\theta)\}\,\pi(\d\theta)=1/2$$
and
$$\int\int (\theta-x)^2\,m(\d x,\d\theta)=\sigma^2,$$
so $w=\sigma^{-2}/2.$
Hence, this calibration idea yields the ``correct" value of $\sigma^{-2}/2$ in this case. This construction requires the user specification of a joint density $m(\d x, \d \theta)$ which 
in some circumstances may prove difficult. One further suggestion is to replace the prior evaluation of the expected 
datum-loss with the observed unit information loss given $x$, 
\begin{equation}\label{eq:euil}
\int\int l(\theta,x) \, m(\d x , \d\theta) \approx \frac{1}{n-p} \sum_{i=1}^n l(\widehat{\theta}_x, x_i) 
\end{equation}
where $$\widehat{\theta}_x = \arg \min_{\theta} \sum_{i=1}^n l(\theta, x_i) $$ is the data-loss estimate of $\theta$ and $p$ is the dimension of $\theta$. For instance, in the above example this leads to, 
$w = \widehat{\sigma}^{-2}/2$
where $$\widehat{\sigma}^2 = \sum_{i=1}^n (x_i - \bar{x})^2/(n-1).$$

It is interesting to note in the above that if it is thought the appropriate choice for $\pi(\theta)$ is flat, possible if the $\theta$ space is bounded, then clearly $\log\{\pi(\widehat{\theta})/\pi(\theta)\}=0$. Thus, to be coherent, we would equally believe $\int l(\theta,x)\,m(\d x|\theta)$ does not depend on $\theta$, where $m(\cdot|\theta)$ is a belief distribution for $x$ given $\theta$. This is a condition which would be hard to justify, as  it would then be also for the uniform prior for $\theta$. If one is used, then we only recommend the value of $w$ is not assigned in the above way.

\vspace{0.2in}
\noindent
{\bf 3.3 Hierarchical loss.} Another way to proceed is to extend the loss function to include $w$ as an unknown parameter. Standard ideas here would suggest we take 
$$L(\theta, w;x,\pi)=w\,l(\theta,x) + \xi l(w) -\log\pi(\theta,w)$$
for some $\xi>0$. We would appear to be making no progress since we now have a $\xi$ to assign. However, this is akin to the hierarchical Bayesian model where uncertainty is propagated via hyper-prior distributions to robustify  the ultimate prior choice at some level. Hence, the allocation of a $\xi$ would not be as crucial as the assignment of a $w$. 

For example, as $w$ is a scale parameter on loss-to-data, taking $l(w)= \log w$  the solution is given by
$$\widehat{\nu}(\theta,w|x,\pi)\propto w^\xi\,\exp\{-w\,l(\theta,x)\}\,\pi(\theta,w)$$
and given that $w^\xi$ can be absorbed in to the prior $\pi$ it is reasonable to assess $\xi$ subjectively. That is, it seems unreasonable to accept that $\pi$ can be chosen subjectively but that $\xi$ can not. 

\vspace{0.2in}
\noindent
{\bf 3.4 Operational characteristics and subjective calibration} The idea here is to set $w$ so that the posterior quantiles are calibrated at some level of error to frequentist confidence intervals based on the estimation of $\theta$ via
minimizing the loss
$$\sum_{i=1}^n l(\theta,x_i).$$
So, if $C_\alpha (w, x_1, \ldots , x_n )$ is the $100(1-\alpha)$\% level confidence interval for $\theta$, then we could select the $w$ such that the posterior distribution of $\theta$, with parameter $w$, is such that
$$\mbox{P}(\theta \in C_\alpha(w,x_1,\ldots,x_n)|x_1,\ldots,x_n) = 1-\alpha.$$
See, for example, the review article by Datta and Sweeting (2005) for references to probability matching priors and posteriors, and Ribatet et al. (2009) for ideas in pseudo-Bayesian approaches with composite likelihoods.

More generally we can consider the subjective setting of $w$ where knowledge of the frequentist sampling statistic of  $\sum_{i=1}^n l(\theta,x_i)$ can assist. To begin note that $w$ is explicitly related to the Bayes Factor quantifying the posterior to prior odds,
$$
\log \left( \frac{\pi(\theta | x)}{\pi(\theta' | x)} /  \frac{\pi(\theta)}{\pi(\theta')} \right) = - w [l(\theta, x) - l(\theta',x)]
$$
where $w [l(\theta, x) - l(\theta',x)]$ measures the update in beliefs in favour of $\theta$ from $\theta'$ on observing $x$. Clearly the larger the difference $[l(\theta, x) - l(\theta',x)]$ the greater the relative evidence in favour of $\theta$, with $w$ determining the scale for unit change. It is interesting to note that should the Bayes Factor be known for any three points $\{\theta, \theta', x\}$ in the joint parameter sample space, $\Omega_{\theta^2} \times \Omega_{X^n}$, then $w$ would be fixed. The idea here is that the analyst is free to contemplate any specific values $\{\theta, \theta', x\}$ for which the distribution of the statistic $S=[l(\theta, x) - l(\theta',x)]$ may be known, and use this knowledge in turn to help elicit a Bayes Factor therefore setting $w$. A concrete example will help:

Suppose $\theta_0$ denotes the unknown mean of a population with prior $\pi(\theta_0) = N(0, v)$ and loss $l(\theta, x) = \sum_i (\theta - x_i)^2$. Consider the design points $\{\theta=\bar{x}, \theta'=0, x\}$ so that the statistic, $S$, is then
$$
S = \sum_i x_i^2 - \sum_i (\bar{x} - x_i)^2
$$
the difference in the sum of squares to the sum of squares around the mean, and log Bayes factor
$$
\log \left( \frac{\pi(\theta | x)}{\pi(\theta' | x)} /  \frac{\pi(\theta)}{\pi(\theta')} \right) = - w S.
$$
The analyst is free to contemplate any value of $n$ and any $x = \{x_1, \ldots, x_n\}$ to help in the elicitation. Let $n$ be chosen large and contemplate $x$ such that the $(1-\alpha) \%$ confidence interval for the unknown mean touches $\theta'=0$. In this case, for large $n$, we know $S = F_{1,n-1}^{-1}(1-\alpha)$, where $F$ denotes the F distribution. If the analyst is prepared to say how their prior beliefs would be updated on observing $x$ in knowledge of this symmetric CI for $\theta_0$ then the $w$ can be set via
$$
w = - \log(BF) / S.
$$
We give a concrete illustration of this approach in Section 5.

\vspace{0.2in}
\noindent
{\bf 3.5 Conjugate loss prior.}
If prior beliefs about $\theta$ can be expressed in the form
$$\pi(\theta)\propto \exp\{-\lambda\,l(\theta,\mu)\}$$
for given parameters $(\lambda,\mu)$, then the posterior has a conjugate type property; that is
$$\pi(\theta|x)\propto \exp\{-wl(\theta,x)-\lambda l(\theta,\mu)\}.$$
Thus the prior has interpretation of prior observation $\mu$ with precision $\lambda$. Thus $\mu$ and $\lambda$ would be standard objects for a Bayesian to specify. If the prior can then be established as the equivalent of $m$ observations, then we obtain $w$ via $w:\lambda=1:m$.   

If the prior is thought not to be able to be specified in such a way, then a good approximation to any prior can be found with choices of $(M,(\mu_j),(\lambda_j))$  such that
$$\pi(\theta)\propto \exp\left\{-\sum_{j=1}^M \lambda_j\,l(\theta,\mu_j)\right\}.$$
If we now write 
$$\pi(\theta|x)\propto 
\exp\left\{-w\, l(\theta,x)-\Lambda\sum_{j=1}^M (\lambda_j/\Lambda)\,l(\theta,\mu_j)\right\},$$
where $\Lambda=\sum_{1\leq j\leq M} \lambda_j$, then we see that now $w:\Lambda=1:m$.

Thus there is an apparent new concept here in that the experimenter is required to think about how  much information, in the form of the number of prior observations, there is available. However, not completely new, since in some conjugate problems there are parameters which do have the interpretation of a prior sample size; the exponential family, for example.

\vspace{0.2in}
\noindent {\bf 4. General forms of information.} In this section we discuss more general forms of information $x$ rather than assume it arises from some unknown $F_0(x)$. The argument is that provided $l(\theta,x)$ has been specified then an update of a belief distribution about $\theta$ is available. Clearly  this does not rely on any assumption about where $x$ came from or indeed how it became known. 

In particular, we provide a definition of conditional probability when non--stochastic information is available. This allows for posteriors to be applied in much more general settings than Bayesian models, which require a stochastic $x$.

%%%%%%%%%%%%%%%%%%%%%%%%%%%%%%%%%%%%%%%%%%%%%%%%%%%%%%%%%%%%%%%%%%%%%%%%%%
%%%%%%%%%%%%%%%%%%%%%%%%%%%%%%%%%%%%%%%%%%%%%%%%%%%%%%%%%%%%%%%%%%%%%%%%%%
%%%%%%%%%%%%%%%%%%%%%%%%%%%%%%%%%%%%%%%%%%%%%%%%%%%%%%%%%%%%%%%%%%%%%%%%%%
%%%%%%%%%%%%%%%%%%%%%%%%%%%%%%%%%%%%%%%%%%%%%%%%%%%%%%%%%%%%%%%%%%%%%%%%%%
%%%%%%%%%%%%%%%%%%%%%%%%%%%%%%%%%%%%%%%%%%%%%%%%%%%%%%%%%%%%%%%%%%%%%%%%%%
%%%%%%%%%%%%%%%%%%%%%%%%%%%%%%%%%%%%%%%%%%%%%%%%%%%%%%%%%%%%%%%%%%%%%%%%%%

\vspace{0.2in}
\noindent
{\bf 4.1 Conditional probability distributions and non-stochastic data.} The theory of conditional probability distributions is a well-established mathematical theory which
provides a procedure to update probabilities taking into account new information.
Such a procedure is available {\sl only} if the information which is used to update the probability concerns stochastic events; that is, events to which a probability is pre-assigned.
In other words, such information needs to be already included into the probability model. In this section, we shall show how the updating approach can be used to define conditional probability distributions based on non--stochastic information.

Information about $\theta$ may arrive in the form of non--stochastic data; such as if an expert declares 
\begin{equation}
I=``\theta\,\,\, \mbox{is close to 0}".\label{nonsinf} 
\end{equation}
This type of information has been discussed by a number of authors and is known to be problematic for the Bayesian especially when such information arises after or during the arrival of stochastic observations $(x_i)$. We cite the paper by Diaconis and Zabell (1982) and in particular refer the reader to example in Section 1.1 of their paper.

We denote by $I$ a piece of information for which no probability model for each $\theta$ is assigned, in other words it is not and can not be determined to be stochastic in any way. To undestand expressions such as this it is worthwhile to recall the type of information we envisage of the type (\ref{nonsinf}).

While a probability model can not connect (\ref{nonsinf}) and $\theta$, they can be connected via a loss function without much difficulty. For example, 
$l(\theta,I)=w\,\theta^2$
for some $w>0$ could be deemed appropriate. Note here we use $I$ to denote information now, replacing the stochastic $x$.
The update $\widehat{\nu}(\theta)$ based on $I$ and $\pi$ can then be considered as a means of defining an operational conditional probability distribution in the presence of non-stochastic information, given by
$$\widehat{\nu}(\theta|I)=\frac{\exp\{-w\,l(\theta,I)\}\,\pi(\theta)}{\int  \exp\{-w\,l(\theta,I)\}\,\pi(\d\theta)}.$$ 
So, the general Bayesian approach introduced in Section 2
provides a general definition of conditional distributions based on non-stochastic information.

For literature on paradoxes related to forcing non-stochastic events into a probability model with a determination of all the alternatives to $I$ we refer the reader to Freund (1965), Gardener (1959), Bar-Hillel and Falk (1982), and Hutchison (1999, 2008).

\vspace{0.2in}
\noindent{\bf 4.2 Partial information.} 
% We now consider a partial information problem. 
% Following from the previous section we label generic information as $I$, whether it be stochastic or non-stochastic. So $I$ could be a statement from an expert relating to the parameter $\theta$ or could be a variable; i.e. $I=X$, assumed to come from some family of distribution functions.

As noted in Section 2, while the parameter of interest is $\theta$, the information $I$ collected may be more informative.  That is, there is within $I$ information which does not assist with the learning about $\theta$, for which it is possible to identify $I_\theta\subset I$ which provides {\sl all} the information about $\theta$. One is therefore interested in constructing the loss function $l(\theta,I_\theta)$, leading to
\begin{equation}
\widehat{\nu}(\d\theta)\propto \exp\{-l(\theta,I_\Theta)\}\,\pi(\d\theta).\label{postpart}
\end{equation}
The partial likelihood, or partial self-information loss, used in proportional hazards model is one such example. While Bayesian practitioners may have adopted such a procedure in the past it would be regarded as lacking motivation. On the other hand, our point is that (\ref{postpart}) represents a valid update of beliefs. We illustrate this approach in Section 5.

\vspace{0.2in} \noindent {\bf 5. Illustrations.} In this section we discuss the application of our approach to important inferential problems. The first problem is one from survival analysis where we have a well motivated proxy likelihood based on partial information. The second example is from model-free clustering where we have a general loss function so that calibration of $w$ is important. A third example, to be found in the Supplementary Material, is for joint inference on a set of quantiles.
In all cases we claim that the choice of loss function is well founded (and unique) and that there is no traditional Bayesian interpretation of the updates we are implementing. Yet the updates we employ do allow us to learn about the specified parameters of interest. All of the models used to generate results are available as open source code in R or Matlab.

\vspace{0.2in}
\noindent
{\bf 5.1 Colon cancer genetic survival analysis.}
Colon cancer is a major worldwide disease with increasing prevalence particularly within western societies. Exploring the genetic contribution to variation in survival times following incidence of the cancer may shed light into the disease eitiology and underlying disease heterogeneity. To this aim collaborators at the Wellcome Trust Centre for Human Genetics, University of Oxford, obtained survival times on 918 cancer patients with germline genotype data at 100,000's of markers genome-wide. For demonstration purposes we only consider one chromosomal previously identified as holding a potential association signal containing 15,608 genotype measurements. The data table $X$ then has $n=918$ rows and $p=15,608$ columns, where $(X)_{ij} \in \{0,1,2\}$ denotes the genotype of the $i$'th individual at the $j$'th marker. Alongside this we have the corresponding $(n \times 2)$ response table of survival times $Y$ with a column of event-times, $y_{i1} \in {\Re}^+$ and a column of indicator variables $y_{i2} \in \{0,1\}$, denoting whether the event is observed or right-censored at $y_{i1}$.

To explore association between genetic variation and time-to-event we employ a loss function derived under proportional hazards, treating the loss to the baseline hazard as a nuisance parameter.
This is based on the Cox proportional hazard (PH) model, one of the most widely used methods in survival analysis since its introduction in Cox (1972). In this log-linear model the hazard rate at time $t$ for an individual with covariate ${\bf x}= \{x_1, \ldots, x_p\}$ is defined as,
$$
h(t | {\bf{x}}) = h_0(t) \exp \left( \sum_{j=1}^p x_j \beta_j \right )
$$
where $h_0(t)$ is a baseline hazard function. In the seminal work of Cox (1972), $h_0(t)$ is treated as a nuisance parameter (or process) that does not enter into the partial-likelihood for estimating the parameters of interest $\bbeta$.

On the other hand, a Bayesian approach to the Cox model necessarily involves the baseline hazard function. There is a limiting argument for the use of the partial likelihood but this is rarely, if at all, used.  Most common is the finite partitioning of the time axis and using a piecewise constant baseline hazard function. Though typically regarded as a nuisance parameter, the Bayesian must specify a full probability model for it. See Ibrahim et al. (2001), Chapter 3, for details, where they note that the proportional hazards model is obtained under a limiting improper prior on the baseline, but it is not known what affect this has on marginal quantities of interest such as marginal model choice probabilities.

Using our construction we can consider only the order of events as partial-information relevant to the regression coefficients, $\bbeta$, via the cumulative loss function,
\begin{equation}
l(\bbeta , {\bf x}) =  \sum_{i=1}^n \log \left ( \frac{\exp\left(\sum_{j=1}^p x_{ij} \beta_j\right)}{\sum_{l\in R_i} \exp\left(\sum_{j=1}^p x_{lj} \beta_j\right)} \right ),\label{eqn:part_lik}
\end{equation}
where $R_i$ denotes the risk set, those individuals not censored or at time $t_i$, and in this way obtain a conditional distribution $\pi(\bbeta | \bx)$.  We assume, $\beta_j \sim N(0, v_j)$
and set $v_j = 0.5$ for our study, reflecting beliefs that associated coefficients will be modest; although we note that one advantage of our approach is that subjective prior information can be integrated into the analysis.

For evidence of effects; i.e. $\beta_j \ne 0$, we can calculate the general Bayes Factor of association at the $j$ th marker as,
$$
BF_j = \frac{ \int_{\beta_j} \exp \left[ -l(\beta_j | \bx_j) \right] \pi(\beta_j ) \d \beta_j}{ \exp \left[ -l(\beta_j = 0 | \bx_j) \right] }
$$
which involves a one-dimensional integral that we calculate via importance sampling.

\begin{figure}
\centerline{
\includegraphics[scale=0.4]{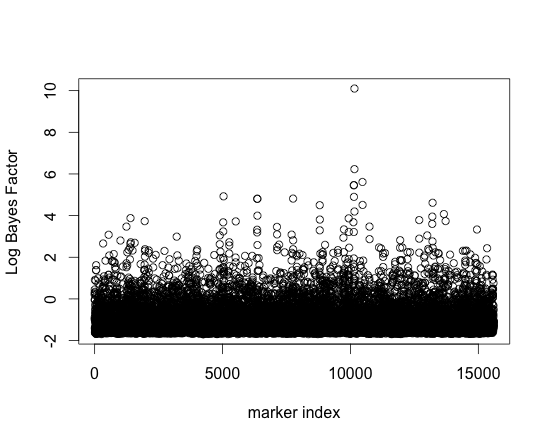}}%{Figs/BF_v_index}}
\caption{Log Bayes Factor vrs marker index along chromosome}\label{BF}
\end{figure}

\begin{figure}
\centerline{
\includegraphics[scale=0.5]{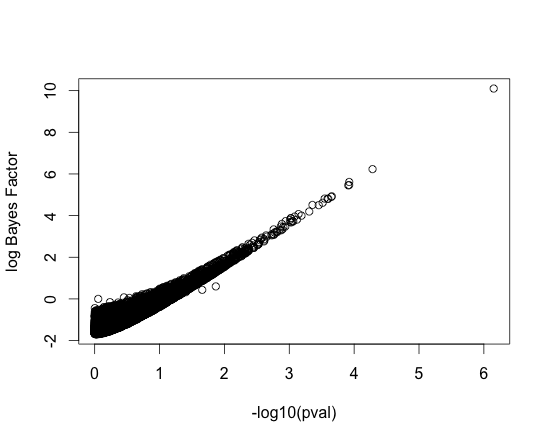}}%{Figs/BF_v_pval_bw}}
\caption{Log Bayes Factor vrs -log10 p-value of association}\label{BFvP}
\end{figure}

\begin{figure}
\centerline{
\includegraphics[scale=0.5]{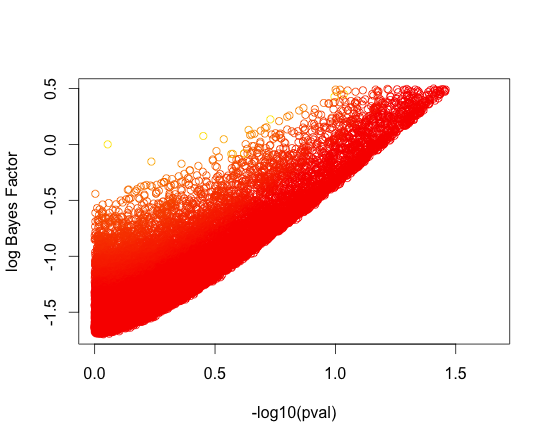}}%{Figs/BF_v_pval_low_col}}
\caption{Log Bayes Factor vrs -log10 p-value of association coloured by standard error in MLE }\label{pvalcol}
\end{figure}

\begin{figure}
\centerline{
\includegraphics[scale=0.5]{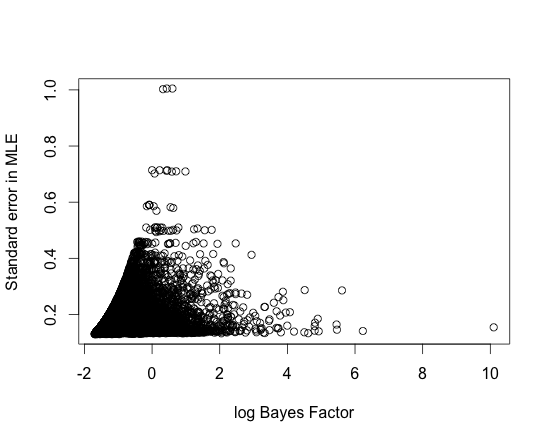}}%{Figs/SE_v_BF}}
\caption{Standard Error in MLE vrs log Bayes Factor }\label{SEvBF}
\end{figure}

We calculated the general Bayes Factors for each marker and in Fig (\ref{BF}) we plot the log Bayes Factors over the chromosome. While there is considerable variation we observe strong evidence of association around marker 10,000. 
% To test if the Laplace approximation is accurate we selected 500 markers at random and ran a Monte Carlo importance sampler with $N(\tilde{\beta}_j, \tilde{\Sigma}_j^{-1})$, and 500 samples. Fig (\ref{LP}) indicates that the Laplace approximation appears accurate. This is not so surprising given we have 918 observations and a single parameter.
%
It is interesting to compare the evidence of association provided by the Bayes Factor Fig (\ref{BF}) in comparison to that obtained using a conventional Cox PH partial-likelihood based test. In Fig (\ref{BFvP}) we plot the log Bayes Factors versus $- \log_{10}$ p-values obtained from a likelihood ratio test. We can see general agreement especially at the markers with strongest association as one would expect for a large sample size. Interestingly there appears to be greater dispersion at markers of weaker association. In Fig (\ref{pvalcol}) we highlight the region of weaker association and colour the points by the standard error of the maximum likelihood estimate. We can see a tendency for markers with less information, greater standard error, to get attenuated towards a logBF of 0 under the general Bayesian approach. This is further highlighted in Fig (\ref{SEvBF}) where we plot the standard error against log Bayes Factors. Markers with high standard error relate to genotypes of rarer alleles and the attenuation reflects a greater degree of uncertainty for association at these markers that contain less information.

\begin{figure}
\centerline{
\includegraphics[scale=0.5]{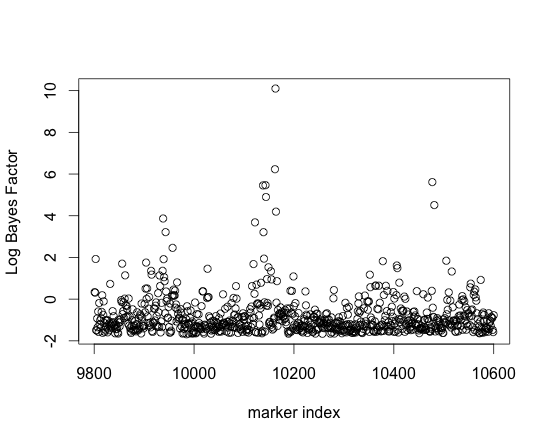}}%{Figs/BF_v_index_hit_region}}
\caption{Log Bayes Factor vrs marker index in the ``hit region'' }\label{BF_hit}
\end{figure}

Returning to the ``hit region'' showing strongest association around marker 10,000, in Fig (\ref{BF_hit}) we see the portion of the graph from Fig (\ref{BF}) containing 800 makers around the marker of strongest association. Due to high colinearity between markers it is not clear whether the signal of association arises from a single effect correlated with others, or from multiple independent association signals. In order to investigate this we developed multiple marker methods. 

We consider a model using potentially all 800 makers in the region and phrase the problem as a variable selection task under a partial-likelihood (loss), in which the user suspects that some of the $p=800$ recorded covariates (\ref{eqn:part_lik}) may not be relevant to variation in survival times.

In the non-Bayesian paradigm, variable selection can proceed by defining a cost function, such as AIC or BIC, that adjusts fit to the data by the number of covariates in the model. Inference proceeds using an optimization algorithm, such as forward or stepwise selection, to find a model that minimises the cost. More recently, penalized-likelihood methods have proved popular (Tibshirani, 1997; Fan and Li, 2002) where the partial-likelihood is maximised subject to some constraint on the norm of the regression coefficients defined by some appropriate sparsity inducing metric.

Despite the enormous impact of Cox PH models and the importance of variable selection, the Bayesian literature in this area is limited.  This is because of the lack of a theoretical foundation to treat $h_0(t)$ as a nuisance parameter, leading to either ad hoc methods or the full specification of a joint probability model. For instance, Faraggi and Simon (1998) and Volinsky et al. (1997) adopt pseudo-Bayesian approaches. The paper of Volinsky et al. (1997) take the BIC as an approximation to the marginal likelihood and they use a branch and bound algorithm to find a set of models with differing sets of covariates with high BIC scores. The difficulty here is that, while the methods are important and well motivated, they are ultimately ad hoc. Moreover, prior information on $\pi({\bm{\beta}})$ does not enter into the calculation of the BIC, meaning that an important aspect of the Bayesian approach is lost.

In contrast, Ibrahim et al. (1999) consider variable selection within a full joint model using a prior specification of a gamma process for the baseline hazard (see also Ibrahim et al. (2001)). This provides a formal Bayesian solution but inference is then conditional on, and sensitive to, the specification of the prior on $h_0(t)$, something the partial-likelihood model explicitly avoids.

\begin{figure}
\centerline{
\includegraphics[scale=0.5]{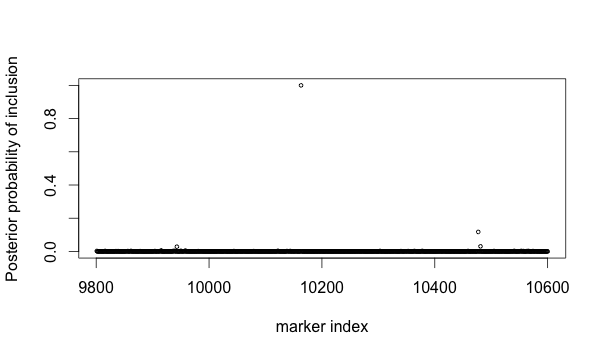}}%{Figs/Post_ssvs_region}}
\caption{Posterior marginal inclusion probability from multiple marker model }\label{Prob_post}
\end{figure}

Here we use the partial-information relevant to the regression coefficients $\bbeta$ via the cumulative loss function (\ref{eqn:part_lik}).
We assume proper priors, $\pi(\bbeta)$ on the regression coefficient,
$$
 \pi(\beta_j) = \left\{
\begin{array}{lll}
  0 &  ~ & {\textrm{if }} \delta_j = 0  \\
  \mbox{N}(0, v_j) & ~ &    {\textrm{otherwise}  },
\end{array}
\right.
$$
where $\delta_j \in \{0,1\}$ is an indicator variable on covariate relevance with,
$\pi(\delta_j) = \mbox{Bin}(a_j)$
and we now treat $\{\delta_1, \ldots, \delta_{800}\}$ as a vector in a joint model. In this way the posterior $\pi(\bdelta | \bx)$ quantifies beliefs about which variables are important to the regression. We use Markov chain Monte Carlo (MCMC) to draw samples approximately from $\pi(\bbeta, \bdelta | \bx)$ from which the marginal distribution on $\bdelta$ can be examined. In particular we make use of an efficient joint updating proposal, $q(\bdelta', \bbeta' | \bdelta)$, within the MCMC as
$
q(\bdelta' , \bbeta' | \bdelta) = q(\bdelta' | \bdelta) q(\bbeta' | \bdelta')
$
where $q(\bdelta' | \bdelta)$ proposes a local move to add, remove, or swap one variable per MCMC iteration in or out of the current model indexed by $\bdelta$, and $q(\bbeta' | \bdelta')$ is a joint independence Metropolis update proposal,
$
q(\bbeta' | \bdelta') = \mbox{N}(\tilde{\bbeta}_{\delta'}, \tilde{{\bm{V}}}_{\delta'})
$
where $\{\tilde{\bbeta}_{\delta'}, \tilde{{\bm{V}}}_{\delta'}\}$ are the MAP and approximate Information Matrix obtained from the combination of log-partial-loss and normal prior. The joint proposal is then accepted with probability,
$$
\alpha  = \min \left\{ 1, \frac{ \exp[-l(\bbeta' | {\bf x})] \pi(\bbeta' | \bdelta') \pi(\bdelta') q(\bbeta, \bdelta | \bdelta') } {  \exp[-l(\bbeta | {\bf x})] \pi(\bbeta | \bdelta) \pi(\bdelta) q(\bbeta', \bdelta' | \bdelta) } \right\}
$$
We ran our MCMC algorithm for 100,000 iterations with prior parameter settings, $\{v_j = 0.5, a_j=1/800\}$, for all $j = 1, \ldots, p$, equivalent to a prior assumption of a single associated marker. In Fig (\ref{Prob_post})  we show the marginal inclusion probability, after discarding 10,000 samples as a burn in. The algorithm showed an overall acceptance rate of 8\% for proposed moves. The model suggest overwhelming evidence for a single marker in the region of index 10200 but also weaker evidence of independent signal in a couple of other regions. 
% We can plot the corresponding log Bayes Factors in Fig (\ref{BF_post}). The additional signals around marker 9900 and 10500 are worthy of follow up and are currently being investigated.
R code to perform the reversible jump MCMC multiple variable sampling for the Cox PH partial-likelihood with normal priors is available on request.

\vspace{0.2in}
\noindent
{\bf 5.2 Bayesian model-free clustering.} Cluster analysis is one of the most widely used and important areas of modern statistics (Hastie et al 2009). In cluster analysis a primary objective is to identify self-similar groups within data, such that observations within a group are deemed more closely related to one another than observations between groups. $K$-means clustering being arguably the most popular clustering method in use today.

The clustering problem is interesting from a formal Bayesian perspective as it raises a number of challenges. The object of interest is the cluster partition mapping, $S$, that allocates observations to clusters. However the partition $S$ as it stands is not a generative model (sampling distribution for observables). To implement clustering the Bayesian analyst is forced to define a sampling distribution for observations within a cluster, $f(x | C_j)$, where $C_j$ denotes parameters associated with the $j$th cluster, with an associated prior probability of cluster membership $p_j$. This leads to the well known marginal mixture representation,
$$
f(x | C) = \sum_{j=1}^K p_j f_j(x | C_j)
$$
the canonical example being with Gaussian mixture components, $f(x | C_j) = N(\mu_j, \Sigma_j)$, which necessitates a further layer of hierarchical priors $\pi(\mu_j, \Sigma_j)$. The mixture model representation to accomplish inference on $S$ has significant consequences. First, cluster membership can be sensitive to the choice of sampling distribution. Secondly Bayesian mixture models suffer from the well known label switching problem (Jasra {\em{et al}} 2005). 

This is a good illustration of where conventional Bayesian updating is restrictive. In cluster analysis the object of interest is typically the partition structure, $\pi(S)$, yet the Bayesian analyst is forced to introduce an additional level of complexity, specifying mixture components with nuisance parameters and priors on nuisance parameters, $\{\mu_j, \Sigma_j\}$, which are of no interest, yet affect the analysis, and define a likelihood which exhibits symmetry to permutation of the labels.

On the other hand, non-Bayesian model-free segmentation methods have a distant advantage in allowing the analyst to concentrate on the object of interest, namely the clustering $S$, typically defined through the specification of a pairwise dissimilarity score between observations $d(x_i, x_j)$. An optimisation algorithm is then used to find the optimal partition $\widehat{S}$ that minimises the score over pairs within and/or between clusters. However, in this approach there is no natural method\footnote{Note, we are not referring to the stability of the optimisation routine to find the global minima but rather that there maybe many competing segmentations of the data with near equivalent scores.} to explore uncertainty, sensitivity or stability in the reporting of $\hat{S}$. Conventional methods of frequentist uncertainty characterisation, via the bootstrap or sub-sampling, cannot be routinely applied to model-free cluster analysis as dissimilarity scores are typically defined over pairs of observations and we only have one realisation of the data. Some methods propose to subsample the data and look at stability of global statistics, such as the choice of the number of groups, but there is generally a lack of theoretical understanding and justification for stability methods in clustering, and moreover such approaches cannot answer the important question of pairwise cluster uncertainty on the joint allocation $P(\{x_i, x_j \in C_k\})$.

In Seldin and Tishby (2001), the authors use PAC Bayesian ideas to consider latent cluster structure to obtain a predictor distribution $Q(y|x)$ for outcome $y$ given $x$ based on a clustering of the $x$. For example, $x$ could be a set of movies, or a set of individuals, or a vector combination of both. The loss function to be minimized is given by
$$L(Q_1,Q_2)=\beta\,\sum_{i=1}^n\int  l(y_i,y)\,Q(d y|x_i)+\int d_{KL}\bigg(Q_2(c|x),Q_2(c)\bigg)\,\lambda(d x)$$
where
$$Q(y|x)=\int Q_1(y|c)\,Q_2(d c|x)$$
and 
$$Q_2(c)=\int Q_2(c|x)\,\lambda(d x)$$
for some probability measure $\lambda$. Here $c$ would denote a clustering of the $x$ and the $(y_i,x_i)$ are the observed data. This minimization is non-trivial and the setting of the coefficient scaling the information in the data can be hard to determine, though the authors suggest a cross-validation approach.  

\begin{figure}
\centerline{
\includegraphics[scale=0.4]{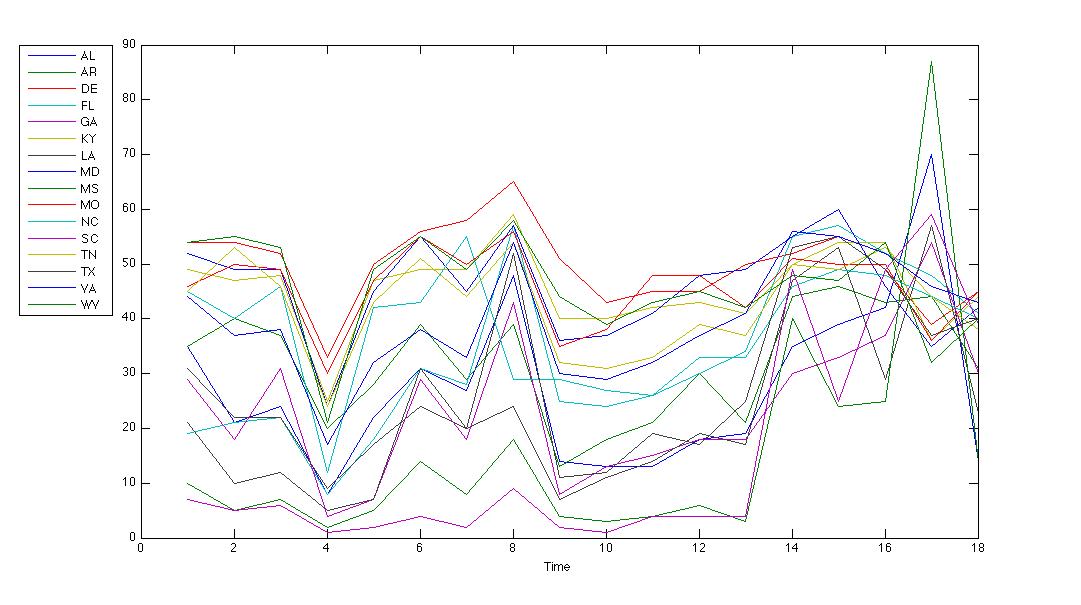}}%{Figs/HartDat2}}
\caption{Voting of Southern States illustrating $\%$ of republican vote for presidential elections every four years beginning in 1900: AL	Alabama; 
AR	Arkansas;
DE	Delaware;
FL	Florida;
GA	Georgia;
KY	Kentucky;
LA	Louisiana;
MD	Maryland;
MS	Mississippi;
MO	Missouri;
NC	North Carolina;
SC	South Carolina;
TN	Tennessee;
TX	Texas;
VA	Virginia;
WV	West Virginia. }\label{fig:vote}
\end{figure}

Our setting is slightly different since we only observe $(x_i)$ and wish to make inference on the underlying clustering of observations into groups. Hence we define a prior distribution directly on the partition, $\pi(S)$, and a loss function $l(S,x_1, \ldots, x_n))$ and use general Bayesian updating. To illustrate this we consider uncertainty analysis of a classic data set considered in Hartigan (1972), illustrated in Figure \ref{fig:vote},  in his highly influential paper that introduced biclustering. Biclustering refers to the simultaneous clustering of observations and covariates (rows and columns) of a data matrix. It having proved extremely useful in modern application areas, particularly in genomics (Cheng \& Church (2000), Tanay {\em{et al}} (2002), Heard {\em{et al}} (2005)). 

Hartigan's paper considered the percentage Republican presidential vote of sixteen southern States in the US over 18 elections covering the years 1900--1968. Hartigan treated the time series as independent covariates in his co-clustering approach. Here, for simple illustration, we maintain the time series ordering, so that the co-clustering is akin to clustering multiple change-point time series with common but unknown change points. We assume the cluster memberships are constant over time, but the time series change at specific breakpoints. Our loss function is defined as in Hartigan (1972) using a sum of squares decomposition,
$$
l(S,x_1, \ldots, x_n) = w \sum_{C \in S} \sum_{i j \in C} (x_{i j} - \bar{x}_C)^2
$$
where $C$ denotes a grouping of States over a particular time period and $\bar{x}_C$ denotes the mean vote of the States and elections assigned to the $C$'th cluster, and posterior distribution 
$
P(S | x) \propto \pi(S) \,\exp\{-w\, l(S, x)\}.
$
The setting of the loss parameter $w$ is a crucial part of the model specification. Following the procedures discussed in Section 3, it's difficult to consider a conjugate specification or a unit information prior on the discrete structures. We instead propose to use a frequentist calibration approach in the following manner. Recall that under a flat prior on $S$ we can set $w$ via a subjective assessment of the posterior ratio at a reference point, 
$$
\frac{P(S | x)}{P(S'|x)} = \exp\{-w [l(S, x) - l(S',x)]\}
$$
and where we can solve for $w$ if all the other elements are given. In elicitation of $w$ we propose to make use of classical results from ANOVA. We take as our first reference point the null partition using a single global cluster, so that the loss $l(S, x) = \sum_i \sum_j (x_{ij}-\bar{x})^2$ is simply the sum of squares around the mean. Then consider a randomised data partition $\{x, S'\}$ that allocates the data uniformly at random to $k$ clusters. Under this allocation scheme we expect,
$$\frac{[l(S,x)-l(S', x)]/k-1}{l(S,x)/(n-k)} \sim F_{k-1,n-k}$$ 
where $F$ denotes the $F$ distribution. We can then use the F-distribution to help in the calibration. For example, if we consider a point in the tails of $F$ such that,
$
f^*_{\alpha} = F^{-1}_{k-1,n-k}(\alpha)
$
with $\alpha \in (0,1)$, and specify
$l(S',x) = l(S,x)/[1 + f^*(k-1)/(n-k)]$
then $l(S',x)$ represents the value of loss such that a randomised allocation has probability $ 1- \alpha$ of producing a smaller loss. Equivalently, with probability $\alpha$ a random allocation would lead to a reduction in loss as high as $[{1 + f^*_{\alpha}(k-1)/(n-k)}]$ relative to the single cluster. When $\alpha$ is large we can be confident that a partition achieving a loss of $l(S',x)$ represents a significant clustering. The analyst can then calibrate $w$ in the following way:
\begin{itemize}
\item Define a reference value for $R=P(S | x)/P(S'|x)$ under a uniform prior, setting $R$ small, say $R=0.001$. 
\item Define a p-value value, $\alpha$, such that should a partition $S'$ achieve a relative reduction in loss of $[{1 + f^*_{\alpha}(k-1)/(n-k)}]$, relative to the global partition $S$, then you would assign relative posterior beliefs of $R$.
\item Set $w = -\log(R) l(S,x)f^*_{\alpha}(k-1)/(n-k)$
\end{itemize}

For the election data we found $w$ is quite stable to the calibration choice of $\{\alpha, R\}$, e.g with $k=3$ we find $w=0.0036$, $\{\alpha=0.99, R=0.01\}$ and $w=0.0012$ for $\{\alpha=0.999, R=0.01\}$. We choose $w=0.0012$ and ran an MCMC algorithm for 100,000 iterations using a burn in of 50,000. The iteration numbers were chosen after experimentation to deliver stable results over multiple runs. The MCMC algorithm was re-run for differing numbers of partitions of States and differing number of time series change points. Table \ref{tab:res} presents the results of the average loss achieved over each run alongside the estimate of the posterior probability for each configuration shown in brackets using a Poisson(3) prior on the number of groups and a Poisson(2) prior on the number of time groupings ($=k_t +1$ in the table), which is the number of change-points +1. Note that the first column in Table \ref{tab:res} equates to standard clustering of the States with zero change points, whereas the first row represents a multivariate change point model. We can see from Table \ref{tab:res} that there is strong evidence for clustering in both time and across States. The maximum posterior probability favours the model with 3 groups of States and 3 time groupings.

\begin{figure}
\centerline{
\includegraphics[scale=0.4]{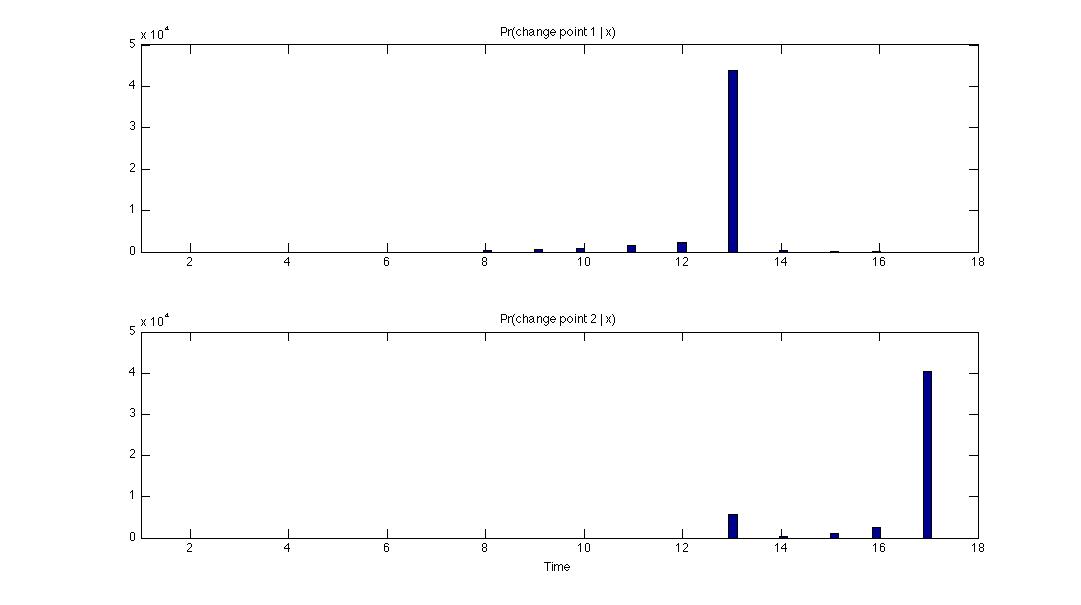}}%{Figs/ChangePoints}}
\caption{Time change point locations for model with 3 groups and 2 change-points }\label{fig:cp}
\end{figure}

\begin{figure}
\centerline{
\includegraphics[scale=0.4]{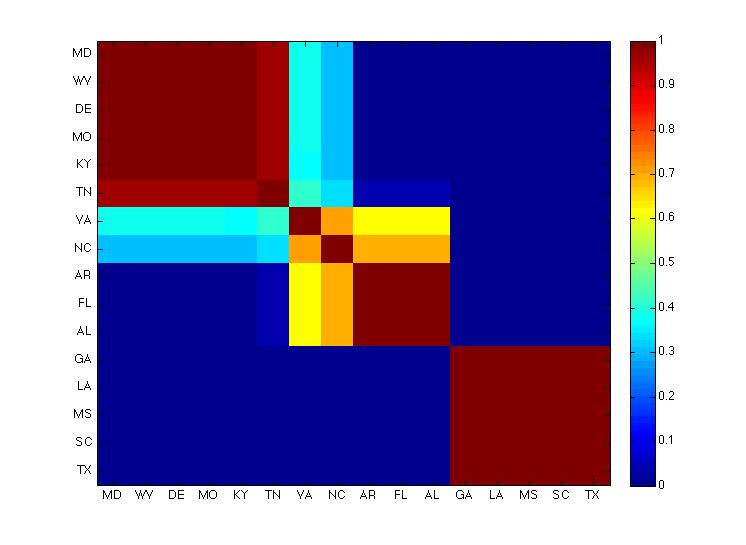}}%{Figs/CoClust}}
\caption{Pairwise co-clustering probabilities across 3 groups and 2 time change-points: see Figure \ref{fig:vote} for labels.}\label{fig:cc}
\end{figure}

We can investigate uncertainty in the partitions and in the cluster allocation of the maximal posterior model. To illustrate this we plot in Figure \ref{fig:cp} the distribution of the location of time-series change points for the $\{k_s=3,k_t=2\}$ model. We can see strong evidence that the change points occur late in the series, which is visually supported by the data in Figure \ref{fig:vote}. The pairwise co-clustering probabilities of this model are shown in Figure \ref{fig:cc}, where each element represents the pairwise probability events $\sum_{S} P(S|x) I[C(x_i)=C(x_j) | S]$, where $C(x_i)$ is the cluster index for the $i$'th State. The cluster blocks show strong concordance with the single co-cluster reported by Hartigan, see Figure 6a in Hartigan (1972). However our method highlights considerable uncertainty in the pairing of Virginia and North Carolina, something we're able to quantify using our general Bayesian approach.
% are two clear blocks of States which exhibit high stability. 
%
%In Seldin and Tishby (2001), the authors use PAC Bayesian ideas to obtain a predictor distribution $Q(y|x)$ for outcome $y$ given $x$ based on a clustering of the $x$. For example, $x$ could be a set of movies, or a set of individuals, or a vector combination of both. The loss function to be minimized is given by
%$$L(Q_1,Q_2)=\beta\,\sum_{i=1}^n\int  l(y_i,y)\,Q(d y|x_i)+\int d_{KL}\bigg(Q_2(c|x),Q_2(c)\bigg)\,\lambda(d x)$$
%where
%$$Q(y|x)=\int Q_1(y|c)\,Q_2(d c|x)$$
%and 
%$$Q_2(c)=\int Q_2(c|x)\,\lambda(d x)$$
%for some probability measure $\lambda$. Here $c$ would denote a clustering of the $x$ and the $(y_i,x_i)$ are the observed data.
%
%This minimization is non-trivial and the $\beta$ hard to determine, though the authors suggest a cross-validation approach. Additionally, the motivation for minimizing $L(Q_1,Q_2)$, since it minimizes an upper bound for a valid but unminimizable loss, is lacking full justification.
%
%Our setting is slightly different since we only observe $(x_i)$ and wish to cluster these into groups. 

\vspace{0.2in} 
\noindent 
{\bf 6. Discussion.} We have provided a basis for general learning and the updating of information using belief probability distributions. Loss functions constructed on spaces of probability measures allow for coherent updating. Specifically, information is connected to the parameter of interest via a loss function and this is the fundamental concept, replacing the restrictive connection based on probability models. We can recover precisely the traditional updating rules such as the Bayes rule when we select the self--information loss function, when it is appropriate to do so.  

The assumptions we make are minimal. That information can be connected to unknown parameters via loss functions and that individuals then act rationally by minimizing their expected  loss. If information is assumed to come from some probability model then we can accommodate  this within our framework by appealing to the self--information loss function equivalent to the negative log-likelihood and so we can argue that loss functions are sufficient for learning mechanisms currently in use.

More generally, we can use loss functions currently employed in a classical context for robust estimation; for example, generalized estimating equations. We can also deal appropriately with partial information where it is only a part of some observed information is useful or relevant for learning about the decision making process based on a particular relevant parameter of interest.

We have developed a rigorous approach to updating beliefs where we are required only to think about which is the best parameter from a chosen model needed to make a decision rather than have to think about a non--existent true model parameter which coincides with the true data generating mechanism.

We believe it is more fundamental to identify parameters of interest through loss functions. The alternative route through a probability model is, we argue, highly restrictive and leads to narrow types of Bayesian updating.  The necessary supporting theory for us is minimal, the construction and minimization of loss functions, whereas for the use of probability models is more intricate and restrictive.

\vspace{0.2in}
\noindent
{\bf Acknowledgements.} The authors are grateful for the detailed comments of two anonymous reviewers and an Associate Editor on a previous version of the paper.

%\end{document}

\vspace{0.2in} \noindent {\bf References.}

\begin{description}

\item Ali, S.M. and Silvey, S.D. (1966). A general class of coefficients of divergence of one distribution from another. {\sl Journal of the Royal Statistical Society, Series B} {\bf 28}, 131--142.

\item Alquier, P. (2008). PAC-Bayesian bounds for randomized empirical risk minimizers. {\sl Mathematical Methods of Statistics}, 17(4), 279-304.

%\item  Allison, P. D. (1995). {\sl Survival Analysis Using the SAS System: A Practical Guide}. Cary NC: SAS %Institute.

%\item Balding, D. J. (2006). A tutorial on statistical methods for population association studies. {\sl Nature Reviews Genetics}, {\bf 7}(10), 781--791.

\item %\bibitem[Bar-Hillel and Falk(1982)]{Bar-Hillel82}
Bar-Hillel, M.~and Falk, R.~(1982).
\newblock Some teasers concerning conditional probabilities.
\newblock \emph{Cognition}, 11:\penalty0 109--122, 1982.

%\item Barron, A., Schervish, M.J and Wasserman, L. (1999)
The consistency of posterior distributions in nonparametric problems
{\sl Annals of  Statistics} {\bf 27}, 536--561.

\item Berger, J.O. (1993). {\sl Statistical Decision Theory and Bayesian Analysis}. Springer Series in Statistics.

%\item Berger, J.O. (2006). The case for objective Bayesian analysis. {\sl Bayesian Analysis} {\bf 1}, 385 -- 402.

%\item Berk, R.H. (1966). Limiting behaviour of posterior distributions when the model is incorrect. {\sl Annals of Mathematical Statistics} {\bf 37}, 51--58.

\item Bernardo, J.M. (1979). Expected information as expected utility. {\sl Annals of Statistics} {\bf 7} 686--690.

\item Bernardo, J.M. and Smith, A.F.M. (1994). {\sl Bayesian Theory}. Wiley.

%\bibitem[Billingsley(1995)]{Billingsley} Patrick Billingsley.

%\item Billingsley, P. (1995).
 %\emph{Probability and measure}.
 %Wiley Series in Probability and Mathematical Statistics. John Wiley
 % \& Sons Inc., New York, third edition, 1995.
 %ISBN 0-471-00710-2.
 %A Wiley-Interscience Publication.

\item Bissiri, P.G. and Walker, S.G. (2010). On Bayesian learning from Bernoulli observations. {\sl Journal of Statistical Planning and Inference} {\bf 140}, 3520-3530.

\item Bissiri, P.G. and Walker, S.G. (2012a). On Bayesian learning from loss functions. {\sl Journal of Statistical Planning and Inference} {\bf 142}, 3167-3173.

\item Bissiri, P.G. and Walker, S.G. (2012b). Converting information into probability measures via the Kullback--Leibler divergence. {\sl Annals of the Institute of Statistical Mathematics} {\bf 64}, 1139-1160.

%\item Box, G.E.P. (1979). Robustness in the strategy of scientific model building, in {\sl Robustness in Statistics}, R.L. Launer and G.N. Wilkinson, Editors. Academic Press: New York.

% \item Box, G.E.P. (1980). Sampling and Bayes' inference in scientific modeling and robustness. {\sl Journal of the Royal Statistical Society, Series A} {\bf 143}, 383--430.

%\item Brown, P.J. and Walker, S.G. (2001). Basyesian priors from loss matching (with discussion). {\sl International Statistical Review} {\bf 80}, 60--82.

%\item Bunke, O. and Milhaud, X. (1998). Asymptotic behaviour of Bayes estimates under possibly incorrect models. {\sl Annals of Statistics} {\bf 26}, 617--644.

%\item Caticha, A. and Giffin, A. (2006). Updating probabilities. In {\sl Bayesian Infer-
%ence and Maximum Entropy Methods in Science and Engineering}, ed.
%by Ali Mohammad-Djafari (ed.), AIP Conf. Proc. 872, {\bf 31}.
%(http://arxiv.org/abs/physics/0608185).

\item Catoni, O. (2003). A PAC Bayesian approach to adaptive classification. Preprint Laborataire de Probabilit\'es et Mod\`eles Al\'eatoires, no. 840.

\item Cesa-Bianchi, N., and Lugosi, G. (2006). {\sl Prediction, learning, and games}. Cambridge University Press.

\item Cox, D. R. (1972). Regression models and life tables (with discussion). {\sl Journal of the Royal Statistical Society, Series B} {\bf 34}, 187-220.

\item Cheng, Y. and Church, G. M. (2000). Biclustering of expression data. {\sl ISMB}  {\bf 8}, pp. 93-103.

\item Datta, G. S., and Sweeting, T. J. (2005). Probability matching priors. {\sl Handbook of statistics}, {\bf 25}, 91-114.

%\item De Blasi, P. and Walker, S.G. (2012). Bayesian asymptotics with misspecified models. {\sl Statistica Sinica} {\bf 23}, 169-187.

%\item Dempster, A.P. (1968). A generalization of Bayesian inference. {\sl Journal of the Royal Statistical Society, Series B} {\bf 30}, 205--247.

\item Diaconis, P. and Zabell, S.L. (1982). Updating subjective probability. {\sl Journal of the American Statistical Association} {\bf 77}, 822--830.

\item Doucet, A., and Shephard, N. (2012). Robust inference on parameters via particle filters and sandwich covariance matrices. {\sl University of Oxford, Department of Economics}. No. 606.

\item de Finetti, B. (1937). La pr\'evision: ses lois logiques, ses sources subjectives. {\sl Annales de l'Institute Henri Poincar\'e} {\bf 7}, 1--68.

%\item Escobar, M.D. (1988). {\sl Estimating the means of several normal populations by nonparametric estimation of the distribution of the means}.
%Unpublished PhD dissertation, Department of Statistics, Yale University.

\item Fan, J and Li, R. (2002). Variable Selection for Cox's proportional Hazards Model and Frailty Model. {\sl Ann. Statist.} {\bf 30}, 1, 74-99.

\item Faraggi, D. and Simon R. (1998). Bayesian variable selection method for censored survival data. {\sl Biometrics} {\bf 54}, 1475-1485.

%\item Feller, W. (1971). %\bibitem[Feller(1971)]{Feller}William Feller.
 %An Introduction to Probability Theory and its Applications.   Vol. II.  Wiley Series in Probability and Mathematical Statistics. John Wiley   \& Sons Inc., New York-London-Sydney, second edition, 1971.

%\item Ferguson, T.S. (1973). A Bayesian analysis of some nonparametric problems. {\sl Annals of Statistics} {\bf 1}, 209--230.

%\item Francesco, P.,  Racugno, W., and Ventura. L. (2011). Bayesian composite marginal likelihoods. {\sl %Statistica Sinica } {\bf 21}, 149--164.

\item %\bibitem[Freund(1965)]{Freund65}
Freund., J.~E. (1965)
Puzzle or paradox?
\emph{Am. Stat.}, 19 (4): 29--44, 1965.

\item %\bibitem[Gardner(1959)]{Gardner59}
Gardner, M.~ (1959).
The Scientific American Book of Mathematical Puzzles and
 Diversions.
 Simon and Schuster, New York, 1959.

%\item Garthwaite, P., Kadane, J. and O'Hagan, A. (2005). Statistical methods for eliciting
%probability distributions. {\sl Journal of the  American Statistical Association} {\bf 100}, 680--701.

%\item Giffin, A. (2008). Maximum Entropy: The Universal Method for Inference. PhD Dissertation,
%University at Albany, State University of New York.

\item Ghosh, J.K. and Ramamoorthi, R.V. (2003). {\sl Bayesian Nonparametrics} Berlin: Springer--Verlag.

\item Goldstein, M. (1981). Revising previsions: A geometric interpretation. {\sl Journal of the Royal Statistical Society, Series B} {\bf 43}, 105--130.

\item Goldstein, M., and Wooff, D. (2007). {\sl Bayes Linear Statistics, Theory \& Methods.} {\bf 716}. Wiley.

%\item Goldstein, M. (2006). Subjective Bayesian analysis: Principles and practice. {\sl Bayesian Analysis} {\bf 1}, 403 -- 420.

\item Hastie, T., Tibshirani, R. and Friedman, J. (2009). {\sl Elements of Statistical Learning}. Springer. 

\item Heard, N. A., Holmes, C. C., Stephens, D. A., Hand, D. J., and Dimopoulos, G. (2005). Bayesian coclustering of Anopheles gene expression time series: study of immune defense response to multiple experimental challenges. {\sl PNAS} {\bf102}(47), 16939--16944.

\item Hirshleifer, J. and Riley, J.G. (1992). {\sl The Analytics of Uncertainty and Information}. Cambridge University Press.

%\item Hjort, N.L., Holmes, C.C., M\"uller, P. and Walker, S.G. (2010). {\sl Bayesian Nonparametrics}. Cambridge University Press.

\item Hoff, P. and Wakefield, J.C.  (2013). Bayesian sandwich posteriors for pseudo-true parameters. {\sl Journal of Statistical Planning and Inference}  {\bf 143}, 1638-1642.

\item H\"uber, P. (1964). Robust estimation of a location parameter. {\sl Annals of Mathematical Statistics} {\bf 35}, 73-101.

\item H\"uber, P. (2009). Robust Statistics (2nd ed.). Hoboken, NJ: John Wiley \& Sons Inc.

\item %\bibitem[Hutchison(1999)]{Hutchison99}
Hutchison, K.~(1999).
\newblock What are conditional probabilities conditional upon?
\newblock \emph{Brit. J. Phi. Sci.}, 50:\penalty0 665--695, 1999.

\item%\bibitem[Hutchison(2008)]{Hutchison08}
Hutchison,K.~(2008).
\newblock Resolving some puzzles of conditional probability.
\newblock \emph{Adv. Sci. Lett.}, 1:\penalty0 212--221, 2008.

\item Ibrahim, J.G. and Chen, M.H. (2000). Power prior distributions for regression models. {\sl Statistical Science} {\bf 15}, 46--60.

\item Ibrahim, J.G., Chen, M.H. and MacEachern, S.N.  (1999). Bayesian variable selection for proportional hazards models. {\sl The Canadian Journal of Statistics}. {\bf 27}, 701-171.

\item Ibrahim, J.G., Chen, M.H. and Sinha, D. (2001). {\sl Bayesian Survival Analysis}, Springer Series in Statistics. 

%\item %\bibitem[Jacod and Protter(2003)]{JacodProtter}
%J.~Jacod and P.~Protter (2003).
%\newblock \emph{Probability essentials}.
%\newblock Springer--Verlag, Berlin, Heidelberg, New York, second edition, 2003.

\item Jasra, A., Holmes, C. C., and Stephens, D. A. (2005). Markov chain Monte Carlo methods and the label switching problem in Bayesian mixture modeling. {\sl Statistical Science}, 50-67.

\item Jiang, W. and Tanner, M.A. (2008). Gibbs posterior for variable selection in high-dimensional classification and data mining. {\sl Ann. Statist.} {\bf 36}, 2207-2231 

\item Kass, R.E. and Wasserman, L.A. (1996).  The selection of prior distributions by formal rules. {\sl Journal of the American Statistical Association} {\bf 91}, 1343--1370.

\item Key, J.T., Pericchi, L.R. and Smith, A.F.M. (1999) Bayesian model choice: What and why? (with discussion). In {\sl Bayesian Statistics 6}, Bernardo, J.M., Berger, J.O., Dawid,
A.P. and Smith, A.F.M. (Eds). Oxford University Press, 343--370.

%\item Kleijn, B.J.K. and van der Vaart, A.W. (2006). Misspecification in infinite dimensional Bayesian statistics. {\sl Annals of Statistics} {\bf 34}, 837--877.

\item Kullback, S. and Leibler, R.A. (1951). On information and sufficiency. {\sl Annals of Mathematical Statistics} {\bf 22}, 79--86.

\item Liang, K.Y. and Zeger, S.L. (1986). Longitudinal data analysis using generalized linear models. {\sl Biometrika} {\bf 73}, 13--22.

\item Marin, J. M., Pudlo, P., Robert, C. P., and Ryder, R. J. (2012). Approximate Bayesian computational methods. {\sl Statistics and Computing}, {\bf 22}(6), 1167-1180.

%\item Lo, A.Y. (1984). On a class of Bayesian nonparametric estimates I. Density estimates. {\sl Annals of Statistics} {\bf 12}, 351--357.

\item McAllester, D. (1998). Some PAC Bayes theorems. In {\sl Proceedings of the 11th Annual Conference on Computational Learning Theory}, ACM, pp.~164--170.

\item Merhav, N. and Feder, M. (1998). Universal prediction. {\sl IEEE Transactions on Information Theory} {\bf 44}, 2124--2147.

\item Muller, U. (2012). Risk of Bayesian inference in misspecified models, and the sandwich covariance matrix. Department of Economics, Princeton University.

\item Ribatet, M., Cooley, D. and Davison, A. C. (2009). Bayesian inference from composite likelihoods, with an application to spatial extremes. {\sl arXiv preprint} :0911.5357.

% \item Rossi, P.H., Berk, R.A. and Lenihan, K.J. (1980). {\sl Money, Work, and Crime: Experimental Evidence}. Academic Press, New York.

\item Royall, R., and Tsou, T. S. (2003). Interpreting statistical evidence by using imperfect models: robust adjusted likelihood functions. {\sl Journal of the Royal Statistical Society: Series B} {\bf 65}, 391--404.

%\item Samuelson, P.A. (1983). {\sl Foundations of Economic Analysis}. (Enlarged Edition). Harvard University Press.

\item Savage, L.J. (1954). {\sl The Foundations of Statistics}. New York. Wiley.

%\item Schepanski, A. and Uecker, W.C. (1984). The value of information in decision making. {\sl Journal of Economic Psychology} {\bf 5}, 177--194.

%\item Shafer, G. (1976). {\sl A Mathematical Theory of Evidence}. Princeton University Press.

\item Seldin, Y. and Tishby, N. (2010). PAC Bayesian analysis of co-clustering and beyond. {\sl Journal of Machine Learning  Research} {\bf 11}, 3595-3646.

\item Shawe-Taylor,  J. and Williamson, R. (1997). A PAC analysis of a Bayesian estimator. In {\sl Proceedings of the 10th Annual Conference on Computational Learning Theory}, ACM, pp.~2--9.

%\item Skilling, J. (1988). The Axioms of Maximum Entropy. In {\sl Maximum-Entropy and
%Bayesian Methods in Science and Engineering}. G. J. Erickson and C. R. Smith
%(eds.) Kluwer, Dordrecht.

%\item Struthers, C.A. and Kalbfleisch, J. D. (1986). Misspecified proportional hazard models. {\sl Biometrika} {\bf 73}, 363--369

\item Tanay, A., Sharan, R., and Shamir, R. (2002). Discovering statistically significant biclusters in gene expression data. {\sl Bioinformatics}, {\bf 18} (1), 136--144.

\item Tibshirani, R. J. (1997). The lasso method for variable selection in the Cox model. {\sl Statistics in Medicine} {\bf 16}, 385-395.

\item Volinsky, C.T., Madigan, D., Raftery, A. E. and Kronmal, A. (1997). Bayesian Model Averaging in Proportional Hazard Models: Assessing the Risk of a Stroke. {\sl Journal of the Royal Statistical Society, Series C}. {\bf 46}, 4, 433--448.

\item von Neumann, J. and Morgenstern, O. (1944). {\sl Theory of Games and Economic Behaviour}. Princeton University Press.

\item Walker, S.G. and Hjort, N.L. (2001). On Bayesian consistency. {\sl Journal of the Royal Statistical Society, Series B}

%\item Walker, S.G. (2004). New approaches to Bayesian consistency. {\sl Annals of Statistics} {\bf 32}, 2028--2043.

%\item White, H. (1982). Maximum likelihood estimation of misspecified models. {\sl Econometrica} {\bf 50}, 1--25.

\item Zellner, A. (1988). Optimal information processing and Bayes's theorem. {\sl The American Statistician} {\bf 42}, 278--284.

%\item Zellner, A. (1997). The Bayesian method of moments (BMOM): theory and applications. {\sl Advances %in econometrics.} {\bf 12}, 85--105.

\item Zhang, T. (2006a). From $\epsilon$-entropy to KL-entropy: Analysis of minimum information complexity density estimation. {\sl Ann. Statist.} {\bf 34}, 2180-2210. 

\item Zhang, T. (2006b). Information theoretical upper and lower bounds for statistical estimation. {\sl IEEE Trans. Inform. Theory} {\bf 52}, 1307-1321. 

\end{description}%\end{thebibliography}

\newpage

\begin{table}[htdp]
\begin{center}
\begin{small}
\begin{tabular}{|c||ccc|}
\hline
%14.4955  -14.3384  -14.7291
%  -13.6921  -13.6482  -13.5786
%  -13.6404  -13.3792  -13.2772
%  -13.9082  -13.5008  -13.4094
 & No. change points in time & $k_t$  & \\ \hline
  & & & \\
 No.  State clusters $k_s$ & $k_t=0$ & $k_t=1$ & $k_t=2$ \\ & & &  \\ \hline
 $k_s=1$ & 7.98 (-14.49) & 6.82 (-14.34) & 6.72 (-14.73)\\ & & &  \\
 $k_s=2$ & 5.36 (-13.69) & 5.13 (-13.65) & 3.19 (-13.58) \\  & & & \\
 $k_s=3$ & 5.09 (-13.64) & 3.92 (-13.38) & 2.36 ({\bf-13.28}) \\ & & & \\
 $k_s=4$ & 4.99 (-13.91)  & 3.32 (-13.50)& 2.02 (-13.41) \\
 \hline
\end{tabular}
\end{small}
\end{center}
\caption{Average loss of partitions $\times 10^4$ across MCMC samples (and log posterior probabilities in brackets). Average loss is $T^{-1}\sum_{i=1}^T l(S_i, x)$ with $S_i \sim \pi(S | x, k_s, k_t)$, where $k_s$ denotes the number of clusters of States and $k_t$ denotes the number of time series change points. Log posterior probabilities {\em{shown in brackets}} using a Poisson(3) and Poisson(2) prior on the number of groups and number of time clusters $=  (k_t+1)$. Maximum posterior clustering shown in bold.}
\label{tab:res}
\end{table}%

%
%\begin{figure}
%\centerline{
%\includegraphics[scale=0.8]{Figs/MC_v_Laplace_500}}
%\caption{Log Bayes Factor using 500 Monte Carlo samples vrs Laplace approximation: at 500 random markers }\label{LP}
%\end{figure}
%

%\begin{figure}
%\centerline{
%\includegraphics[scale=0.8]{Figs/BF_region}}
%\caption{Log Bayes Factor for inclusion in multiple marker model }\label{BF_post}
%\end{figure}

%
%\begin{figure}
%\centerline{
%\includegraphics[scale=0.8]{Figs/rq_v_bq}}
%\caption{Estimates of standard deviation of Bayesian posterior samples versus those for standard errors of the quartiles, estimated from standard Exponential data (red) and Normal (black). }\label{rq}
%\end{figure}

\end{document}